%% file: sstable.tex
\documentclass[11pt]{article}
\usepackage{amsmath}
\usepackage{amsfonts}
\usepackage{amssymb}
\usepackage[dvips]{graphicx}
\usepackage{texdraw}

\newcommand{\C}{\mathbb{C}}
\newcommand{\CP}{\mathbb{CP}}
\newcommand{\cell}{\triangle}
\newcommand{\cod}{\delta}
\newcommand{\cone}{\alpha}
\newcommand{\cross}{\times}
\newcommand{\dsum}{\oplus}
\newcommand{\eopf}{$\square$} 
\newcommand{\equ}{\cong}

\newcommand{\f}{\mathbb{F}}
\newcommand{\face}{\sigma} 
\newcommand{\faceof}{\prec} 
\newcommand{\fan}{\Sigma} 
\newcommand{\gives}{\Rightarrow}
\newcommand{\halfsp}{\textrm{\bf H}}
\newcommand{\inter}{\cap}
\newcommand{\into}{\hookrightarrow}
\newcommand{\iso}{\simeq}
\newcommand{\nfaceof}{\nprec}
\newcommand{\p}{\mathbb{P}}
\newcommand{\polytp}{\cell}

\newcommand{\Q}{\mathbb{Q}}
\newcommand{\R}{\mathbb{R}}
\newcommand{\Rplus}{\R^{+}}
\newcommand{\sheafo}{\sheaf{O}}  
\newcommand{\sminus}{\setminus}
\newcommand{\Span}{\textrm{Span}}
\newcommand{\Spec}{\textrm{Spec}}

\newcommand{\subspace}{\subset}
\newcommand{\subsp}{\subspace}

\newcommand{\tensor}{\otimes}
\newcommand{\triv}{\mathbf{1}}
\newcommand{\union}{\cup}
\newcommand{\Z}{\mathbb{Z}}

\newcommand{\blowup}[1]{\widetilde{#1}}
\newcommand{\bundle}[1]{\mathcal{#1}}
\newcommand{\cl}[1]{\overline{#1}}

\newcommand{\dualcn}[1]{{#1}^{\vee}}

\newcommand{\family}[1]{\mathcal{#1}}
\newcommand{\ideal}[1]{\mathcal{#1}}
\newcommand{\inverse}[1]{{#1}^{-1}}
\newcommand{\linearsys}[1]{|{#1}|}  

\newcommand{\numin}[1]{| {#1} |}
\newcommand{\oomit}[1]{\hat{#1}}
\newcommand{\pb}[1]{\pullback{#1}}
\newcommand{\proof}[1]{{\rm {\bf Proof:} #1 \eopf}}
\newcommand{\pullback}[1]{{#1}^*}
\newcommand{\restricto}[1]{|_{#1}}
\newcommand{\sheaf}[1]{\mathcal{#1}}
\newcommand{\transpos}[1]{{#1}^{t}}
\newcommand{\coho}[2]{H^{{#1}}(#2)}
\newcommand{\combn}[2]{C_{#1}^{#2}}  

\newcommand{\hodgenum}[2]{h^{{#1}, {#2}}}

\newcommand{\pairing}[2]{\langle {#1}, \, {#2} \rangle}
\newcommand{\draw}[1]{\begin{texdraw} \drawdim{cm} #1 \end{texdraw}}
\newcommand{\includeps}[2]{\begin{figure}[!h!tbp]\begin{center}\includegraphics[#1]{#2}\end{center}\end{figure}}
\newcommand{\equlist}[2]{\left. \begin{array}{#1} #2 \end{array} \right.}
\newcommand{\equlabel}[2]{\begin{equation}\label{eq:#1} #2 \end{equation}}
\newcommand{\citeequ}[1]{(\ref{eq:#1})}
\newcommand{\bibli}[2]{\begin{thebibliography}{#1} #2 \end{thebibliography}}
\newcommand{\bitem}[4]{\bibitem{#1} {#2}: \emph{#3}, #4.}
\newtheorem{Lemma}{Lemma}[section]
\newtheorem{Examp}{Example}[section]
\newtheorem{Clm}[Lemma]{Claim}
\newtheorem{Corl}{Corollary}[Lemma]
\newtheorem{Def}[Lemma]{Definition}
\newtheorem{Prop}[Lemma]{Proposition}
\newtheorem{Thm}[Lemma]{Theorem}
\newtheorem{Rmk}{Remark}[Lemma]
\newcommand{\lemma}[2]{\begin{Lemma} \label{Lemma:#1} #2 \end{Lemma}}
\newcommand{\exmp}[2]{\begin{Examp} \label{Examp:#1} {\rm #2} \end{Examp}}

\newcommand{\defin}[2]{\begin{Def} \label{Def:#1} #2 \end{Def}}
\newcommand{\prop}[2]{\begin{Prop} \label{Prop:#1} #2 \end{Prop}}
\newcommand{\thm}[2]{\begin{Thm} \label{Thm:#1} #2 \end{Thm}}
\newcommand{\rmk}[2]{\begin{Rmk} \label{Rmk:#1} \rm{#2} \end{Rmk}}
\newcommand{\citelem}[1]{lemma~\ref{Lemma:#1}}
\newcommand{\citexmp}[1]{example~\ref{Examp:#1}}

\newcommand{\citedef}[1]{definition~\ref{Def:#1}}
\newcommand{\citeprop}[1]{proposition~\ref{Prop:#1}}
\newcommand{\citethm}[1]{theorem~\ref{Thm:#1}}
\newcommand{\citermk}[1]{remark~\ref{Rmk:#1}}

\newcommand{\Citethm}[1]{Theorem~\ref{Thm:#1}}

\newcommand{\concav}{C}
\newcommand{\ctori}{{\C^{*}}}

\newcommand{\flag}{\mathcal{F}}

\newcommand{\simplex}{\kappa}

\newcommand{\upperhalf}{\halfsp^+}

\newcommand{\lat}[1]{{\bf #1}}
\newcommand{\new}[1]{\widetilde{{#1}}}
\newcommand{\support}[1]{| {#1} |}

\newcommand{\shortexact}[7]{{{#1}} \longrightarrow {{#2}} \stackrel{#3}{\longrightarrow} {{#4}} \stackrel{#5}{\longrightarrow} {{#6}} \longrightarrow {{#7}}}
\begin{document}
\title{Semi-Stable Degeneration of Toric Varieties and Their Hypersurfaces}
\author{Shengda Hu\\
        University of Wisconsin-Madison\\
        Email:hu@math.wisc.edu}
\date{}
\maketitle

\begin{abstract}
We provide a construction of examples of semistable degeneration via toric geometry. The applications include a higher dimensional generalization of classical degeneration of $K3$ surface into $4$ rational components, an algebraic geometric version of decomposing $K3$ as the fiber sum of two $E(1)$'s as well as it's higher dimensional generalizations, and many other new examples.
\end{abstract}

\section{Introduction}
Surgery operation has been a central topic in every branch of topology and geometry. For example, in 3-dimension, geometrization conjecture states that one can classify 3-dimension into eight different ``geometries'' up to certain connected sum surgeries. For 4-manifold, we don't have such a classification conjecture precisely because we do not yet understand the surgeries of 4-manifolds. Unfortunately, connected sum is not a surgery in symplectic geometry. Instead, surgery theory in symplectic geometry is more subtle and influenced more by surgery in algebraic geometry. The latter is the main topic of this article.

There are primarily two types of surgeries in algebraic geometry. One class is birational transformation such as blow up-blow down, flip-flop. The birational geometry is  a central topic in algebraic geometry. The famous Mori theory was designed to study this type of surgeries. Their importance in symplectic geometry was documented by many important papers \cite{mcduff}, \cite{liruan}.  Another class of surgery comes from deforming the defining equation of algebraic varieties. Suppose that $F_t=\{f_t=0\}$ for $t\neq 0$ defines a family of smooth algebraic varieties. We can let $t\rightarrow 0$ such that $F_0=\{f_0=0\}$ acquires singularities. We further require the singularities of $F_0$ to be normal crossing. In this case, $F_0$ can be viewed as a union of smooth manifolds intersecting along (complex) codimension 1 submanifold. Namely, we ``decompose'' $F_t$ as a union of manifolds intersecting each other transversely. This is precisely what we mean surgery in topology. However, this is still too general to have any nice theory. Mumford proposed the notion of ``semi-stable degeneration''.

\defin{sstable}{A \emph{degeneration} is a morphism $\pi : X \to U \stackrel{open}{\subset} \C$, and $0 \in U$, such that the restriction $\pi : X \sminus \inverse{\pi}(0) \to U \sminus \{0\}$ is smooth. $\pi$ is \emph{semi-stable}, if $X$ is non-singular, and the fiber $\inverse{\pi}(0)$ is reduced, with non-singular components crossing normally (cf. \cite{kkms} Ch II).}

In many ways, semi-stable degeneration is the correct set-up to work with. For example, the monodromy is unipotent. Well-known Clemens-Schmidt exact sequence related cohomology of $F_t$ to that of $F_0$. Hodge structure has a nice limit at $F_0$. There is a theory  of log differential forms on $F_0$. The list goes on.

In a different direction, semi-stable degeneration also appears in quantum cohomology. During last several years, a central topic in quantum cohomology is the surgery theory for Gromov-Witten invariants. In symplectic category, there are two different versions of surgery theory by Li-Ruan \cite{liruan} and Ionel-Parker \cite{ionelparker}. Recently, an algebro-geometric version analogous to that of Li-Ruan was obtained by J. Li \cite{lijun}. These surgery theories have many important applications in quantum cohomology. The set-up in these surgery theories is precisely semi-stable degeneration. One of purposes of this article is to provide many more interesting examples to apply these surgery techniques, which the author hope to investigate in the future.

Although we have many beautiful theories for semi-stable degeneration, surprisingly, there are few concrete examples. Mumford proved a general theorem that every degeneration induces a semi-stable degeneration. However, Mumford's theorem is not constructive and not very useful in terms of constructing example. The only nontrivial example in author's knowledge is the semi-stable degeneration of $K3$-surface into four copies of blow up of $\CP^2$ along 6-points. It is  a beautiful construction and let us briefly describe it. We first take a degeneration of smooth degree 4 hypersurfaces of $\CP^3$ to a union of four hyperplanes, say a pencil. This pencil has base loci. We can blow up $\CP^3$ to resolve base loci. In this process, we increase the number of components in central fiber. Now the total space is a 3-fold and Mori theory applies. We use Mori theory to blow down the total space. The resulting one is a semi-stable degeneration where the central fibers are four copies of the blow-up of $\CP^2$ along 6-points. In higher dimension, we do not have Mori theory. Hence, we do not know if there is a semi-stable degeneration of quintic Calabi-Yau 3-folds into five copies of blow-up of $\CP^3$. Our method in this article will provide such an example and much more. In fact, we will describe a general method to construct semi-stable degeneration of toric varieties and their hypersurfaces. Then, we will use our methods to construct many examples. Our examples include the classical example of $K3$ and its natural generalizations into higher dimension such as Calabi-Yau hypersurface of $\CP^n$. Many other examples are also constructed.

Toric varieties are the testing ground for many ideas, because they are the next simplest varieties after projective spaces, while at the same time they are not at all trivial. They contain many interesting varieties as subvariety. For example, most Calabi-Yau 3-folds known up to now are subvarieties of toric varieties. When considering degenerations, it's natural to consider toric geometry, since it provides the local picture of a semi-stable degeneration. In fact, it's how toric varieties arise in \cite{kkms}. In this paper, we combine these two views, and use toric varieties as ambient spaces of both the fiber and the total space of a semi-stable degeneration. Toric varieties can be thought of as compactifications of an algebraic torus $(\C^*)^n$, and polytopes or fans are used to describe how the compactifications work. In principle, one could obtain everything about a toric variety from the fan or polytope that defines it. We'll consider projective toric varieties which are those that can be defined by a polytope. In this case, the polytope not only specifies the toric variety, it also specifies a polarization of it, i.e. an embedding into some projective space. This will enable us to write down how the equations change when the toric variety (or its subvarieties) degenerates. On the technical side, fans make it easier to describe morphisms between toric varieties, while polytopes are much better at specifying the data we need to construct the degenerations. Through polytopes, our construction connects back to symplectic geometry via the moment map --- it can be thought of as ``doing surgery with respect to the image of the moment map''.

The first main result in this paper is that we can talk about semi-stable degeneration in toric category. Namely, the total space of the degeneration is a toric variety, the smooth fibers are toric varieties, and the components in the central fiber are toric varieties as well. The initial data is a partition $\Gamma$ of a non-singular polytope $\polytp$ (thus $X$, the toric variety defined by $\polytp$, is non-singular), and the theorem is

\vspace{0.1in}\noindent
{\bf \Citethm{degeneration}} \emph{If $\Gamma$ is a non-singular semi-stable partition, then there exists a toric variety $\new{X}$ and a map $p : \new{X} \to \C$, such that $p$ is a semi-stable degeneration of $X$ to $\inverse{p}(0)$. The dual graph of the fiber $\inverse{p}(0)$ is given by the dual graph of the partition $\Gamma$, and a component in $\inverse{p}(0)$ is the toric variety defined by the corresponding subpolytope in $\Gamma$.}

\vspace{0.1in}
Given a semi-stable degeneration $p : \new{X} \to \C$, we can consider a smooth subvariety $\new{V}$ in $\new{X}$. Suppose that $\new{V}$ intersects transversely with general fiber of $p$ as well as all the components in the fiber $\inverse{p}(0)$. Then in a small neighbourhood of $0$, the restriction of $p$ to $\new{V}$ gives a semi-stable degeneration. In our case, we also see that the components in $\inverse{p\restricto{\new{V}}}(0)$ are subvarieties of toric varieties. One application is

\vspace{0.1in}\noindent
{\bf \Citethm{CYdegen}} \emph{There is a semi-stable degeneration of degree $n+1$ Calabi-Yau hypersurface in $\CP^n$, such that all the components in the singular fiber are rational, and the dual graph of the central fiber is the triangulation of $S^{n-1}$ by the boundary of $n$-simplex.}

\vspace{0.1in}\noindent
Another result also concerning Calabi-Yau hypersurface is a generalization of realizing $K3$ surface as fiber-sum of two $E(1)$'s.

\vspace{0.1in}\noindent
{\bf \Citethm{degenCYodd}} \emph{There exists a semi-stable degeneration of Calabi-Yau hypersurface in $\CP^{2k+1}$ into two components $Y_i$ $(i = 1, 2)$, such that $Y_1 \iso Y_2 \iso Y$ satisfy the following: $(K_Y)^2 = 0$, $-K_Y$ is effective and defines a Calabi-Yau subvariety $Z$ of $Y$, $Z$ is isomorphic to a hypersurface in $\CP^k \cross \CP^k$, and the singular loci of the central fiber is $Z$.}


\vspace{0.1in}
The ideal of the proof of \citethm{degeneration} is the following. Starting from a non-singular semi-stable partition (\citedef{partition}) of $\polytp$, we can construct a concave function $f$ on $\polytp$, which is integral affine on each sub-polytope of the partition. We take $\new{\polytp}$ to be the convex polytope that is cut out from the cylinder $\polytp \cross \R$ by the graph of $f$. If we can choose $\new{\polytp}$ to be non-singular, then the total space of the degeneration is the toric variety $\new{X}$ defined by $\new{\polytp}$. The rest of \citethm{degeneration} follows from toric arguments.

The paper is organized as following: Section $2$ gives a brief introduction to toric varieties. We state results in toric geometry that will be used both in the proof of the main results, and in the examples. Many results in toric geometry has statements in terms of fans, and we'll state them in terms of either fans or polytopes depending on how we're going to use them. Section $3$ is purely combinatorical and it contains the construction of the polytope $\new{\polytp}$ of one higher dimension mentioned above from a non-singular (or mildly singular) partition of $\polytp$. In section $4$, we use the result in section $3$ (\citethm{construction}) to prove \citethm{degeneration} and a similar result on weak semi-stable degeneration of simplicial toric varieties. In the same section (\citermk{degen3}), we also give explicit description in coordinates of the degeneration of toric varieties as subvarieties in some big projective space. This description leads to one of the basic constructions described in $5.1$, which is again used extensively in the examples. Section $5$ contains the examples, including the example on Calabi-Yau hypersurfaces mentioned above ($5.4$), an algebraic geometrical version of $E(1)$ fiber-sum with another copy of $E(1)$ produces a $K3$ surface ($5.3$), and degenerations of degree $d$ hypersurfaces into a chain of $d$ rational varieties ($5.5$). Section $6$ contains a generalization of the construction to give semi-stable degenerations with higher dimensional base.

After the completion of this paper, Mark Gross pointed out to the author that a construction similar to that of section 3,4 are used by V. Alexeev in \cite{alexeev}. In \cite{alexeev}, Alexeev has completely different purpose of constructing the completion of an arbitrary one-parameter family of stable toric pair. In particular, semi-stability was not an issue in \cite{alexeev}. On the other hand, the main purpose of our work is to construct and identify special types of degeneration which are semi-stable. Therefore, we have completely different goals in the mind while using a similar construction.

\vspace{0.1in}
{\bf Acknowledgments}

I want to thank my advisor, Professor Yongbin Ruan, for his guidance and encouragement. I benefitted a lot from discussions with Professor Miles Reid, Rajesh Kulkarni, Wojciech Wieczorek, Mainak Poddar and Bohui Chen. I want to thank them all here.

\section{Brief introduction to toric varieties}

We outline the main results and definitions we'll need in toric geometry and refer to \cite{cox}, \cite{fulton} and \cite{oda} for more detailed statements and proofs.

\defin{toric}{\emph{\cite{fulton} (p.3)} A \emph{toric variety} is a normal variety $X$ that contains a complex torus $(\ctori)^n$ as dense open subset, together with an action $(\ctori)^n \cross X \to X$ of $(\ctori)^n$ on $X$ that extends the natual action of $(\ctori)^n$ on itself.}
While we have this definition, the construction of concrete toric varieties, and proof of many of their properties start from lattices, fans and polytopes.

$M$ and $N$ are dual lattices of rank $n$, $M_{\R} = M \tensor \R$. A \emph{rational convex polyhedra cone} $\cone$ in $N_{\R}$ is a convex cone generated by finite many vectors in $N$. It can be given by:
$$\cone = \{u \in N_{\R} ~ | ~ \pairing{m_i}{u} \geq 0, m_i \in M ~\textrm{ for }~ i = 1, \ldots, s\}.$$
We'll take all the cones in the following to be \emph{strongly convex} meaning that they do not contain any non-trivial linear subspace through origin. We'll refer to such cone simply as \emph{cone} since no other cones will appear in this paper. For a cone $\cone$ in $N_{\R}$, we define the \emph{dual cone} $\dualcn{\cone}$ to be a cone in $M_{\R}$ given by
$$\dualcn{\cone} = \{v \in M_{\R} \, | \, \pairing{v}{u} \geq 0 ~\textrm{ for all }~ u \in \cone\}.$$
A \emph{face} $\tau$ of cone $\cone$ is the intersection of $\cone$ with one of its supporting hyperplanes. A \emph{fan} $\fan$ in $N_{\R}$ is defined to be a collection of cones such that
\begin{enumerate}
\item Each face of a cone in $\fan$ is also a cone in $\fan$,
\item The intersection of two cones in $\fan$ is a face of each.
\end{enumerate}

Classical construction of toric varieties starts with a fan $\fan$ in $N_{\R} = N \tensor \R$. The toric variety $X = X_{\fan}$ is defined by gluing together the affine toric varieties 
$X_{\cone} = \Spec(\C[\dualcn{\cone} \inter M]).$
 Let 
$\phi : N' \to N$
be a homomorphism of lattices, and $\fan'$ is a fan in $N'$. The map is said to be \emph{compatible with the fans} if for each cone $\cone'$ in $\fan'$, there is some cone $\cone$ in $\fan$ such that $\phi(\cone') \subset \cone$ (cf. \cite{fulton} p.22). Such map induces a morphism
$\phi_* : X_{\fan'} \to X_{\fan}$ of toric varieties
which preserves the correspondence between subvarieties and cones, i.e., 
\equlabel{derMorph}{\phi(\cone') \subset \cone ~~ \gives ~~ \phi_* : Y_{\cone'} \to Y_{\cone}}
Let $\support{\fan}$ be the support of $\fan$, i.e. the union of the cones in $\fan$, then (cf \cite{fulton} p.39)
$$\phi_* ~\textrm{ is proper }~ \iff \inverse{\phi}(\support{\fan}) = \support{\fan'}.$$
If $\phi$ is an isomorphism between lattices, and $\fan' = \phi(\fan)$, then $\phi_*$ is an isomorphism between the varieties, which induces a change of coordinates on the embedded big torus. If $\support{\fan} = N_{\R}$, it's called \emph{complete}, and the corresponding $X_{\fan}$ is a compact variety.

The construction can also start with a polytope (not necessarily compact) $\polytp$ in $M_{\R}$. A (compact) \emph{rational convex polytope} $\polytp$ in $M_{\R}$ is the convex hull of finite number of points of $M$. An alternative definition encompassing non-compact polytopes is the following:
$$\polytp = \{v \in M_{\R} ~ | ~ \pairing{v}{n_i} \geq -a_i, ~ n_i\in N, ~a_i\in\Z ~\textrm{ for }~ i = 1, \ldots, s\}.$$
A \emph{face} $\face$ of $\polytp$ is the intersection of $\polytp$ with one of its supporting hyperplanes, and we'll denote it by $\face \faceof \polytp$. For each $\face \faceof \polytp$, the cone $\cone_{\face}$ \emph{dual to} $\face$ is defined by 
$$\cone_{\face} = \{u \in N_{\R} \, | \, \pairing{v}{u} \leq \pairing{v'}{u} ~\textrm{ for all }~ v \in \face ~\textrm{ and }~ v' \in \polytp\}.$$
The dual cone $\dualcn{\cone_{\face}}$ in $M_{\R}$ is generated by vectors $v' - v$ where $v$ and $v'$ vary among the vertices of $\face$ and $\polytp$ respectively. It can be shown that $\cone_{\face}$'s form a fan $\fan_{\polytp}$, which gives a toric variety $X = X_{\polytp}$  (cf \cite{fulton} p.26). Let $\psi : M_{\R} \to M_{\R}$ be an affine transformation which preserves $M$, and $\polytp' = \psi(\polytp)$, then $\fan_{\polytp'}$ is the image of $\fan_{\polytp}$ under an isomorphism $\psi^* : N_{\R} \to N_{\R}$. The varieties $X_{\polytp'}$ and $X_{\polytp}$ are isomorphic. We make the following definitions:
\defin{polytp}{A polytope $\polytp$ in $M_{\R}$ is \emph{simplicial}, if there are exactly $n$ edges at each vertex, and the primary vectors at each vertex span $M_{\R}$ as vector space. A fan $\fan$ in $N_{\R}$ is \emph{simplicial}, if all the cones in $\fan_{\polytp}$ is simplicial. A vertex $p$ of a simplicial $\polytp$ is \emph{non-singular} if the primary vectors at $p$ span the lattice $M$. $\polytp$ is a \emph{non-singular} polytope (not necessarily compact), if $\polytp$ is simplicial and all its vertices are non-singular.}
We then have the following:
\prop{toricFacts}{
\begin{enumerate}
\item $\polytp$ is compact $\gives \fan_{\polytp}$ is a complete fan $\gives X_{\polytp}$ is compact, 
\item $\polytp$ is simplicial $\gives \fan_{\polytp}$ is simplicial $\gives X_{\polytp}$ is an orbifold,
\item $\polytp$ is non-singular $\iff X_{\polytp}$ is non-singular.
\end{enumerate}}
A toric variety can be constructed from a polytope in another way (cf. \cite{cox}). Let $x = (x_1, \ldots, x_n) \in (\ctori)^n$, and $x^{\lat{m}} = \prod_{i = 1}^n x_i^{m_i}$ where $\lat{m} = (m_1, \ldots, m_n) \in M$ is a lattice point. Define an embedding of $(\ctori)^n$ into $\CP^{\ell - 1}$ by
$$i_{\polytp} : x \mapsto [x^{\lat{m}_1} : \ldots : x^{\lat{m}_{\ell}}],$$
where $\polytp \inter M = \{\lat{m}_1, \ldots, \lat{m}_{\ell}\}$. The closure of the image of $i_{\polytp}$ is a toric variety $X$. An example of this construction is the Varonese embedding of $\CP^m$ into higher dimensional projective spaces.

A special family of subvarieties of a toric variety $X$ is given by the closure of the $(\ctori)^n$ orbits, which are one-to-one correspond to the cones in the fan $\fan$ that defines $X$, and they are called \emph{invariant subvarieties}. The correspondence reverses the inclusions. Let $V(\cone)$ denote the invariant subvariety corresponding to $\cone$. Let $A_k(X)$ be the Chow group of $X$, then (cf \cite{fulton} p.96)
\prop{chowgp}{$A_k(X)$ is generated by the classes of the invariant subvarieties that corresponds to $n-k$ dimensional cones in $\fan$.}

Suppose that $X$ is given by a compact simplicial polytope $\polytp$, then $X$ is projective. By passing to the fan $\fan_{\polytp}$, we have a correspondence of invariant subvarieties with the faces of $\polytp$ which respects inclusion. Let $V(\tau)$ ($\cone_{\tau}$) be the invariant subvariety (cone in $\fan_{\polytp}$) corresponds to $\tau \faceof \polytp$. We have (cf. \cite{fulton} \S5.1 \S5.2)
\prop{cohoRing}{$A^k(X)_{\Q} = A_{n-k}(X) \tensor \Q \iso H_{2(n-k)}(X, \Q) \iso H^{2k}(X, \Q)$, and the intersection of cycles makes the Chow groups into a graded ring. For different faces $\sigma$ and $\tau$ that don't have each other as faces, $V(\sigma) \cdot V(\tau) \neq 0 \iff \sigma \inter \tau \neq \emptyset$.}
When $X$ is nonsingular, we can take the coefficient group to be $\Z$, and the intersection takes the form
$$V(\sigma) \cdot V(\tau) = V(\gamma) \iff \sigma \inter \tau = \gamma.$$

The divisor on $X$ corresponding to a $(n-1)$-face $\rho \faceof \polytp$ will be denoted $D_{\rho}$. Let $D = \sum_{\rho} a_{\rho} D_{\rho}$ be a Cartier divisor on $X$. Then $D$ defines a piecewise linear function $\phi_D$ on $N_{\R}$ by setting $\phi_D(v_{\rho}) = -a_{\rho}$ where $v_{\rho}$ is the primary vector of the ray corresponding to $\rho$, then extend so that $\phi_D$ is linear on each cone of $\fan_{\polytp}$.
\defin{supptFun}{$\phi_D$ as defined above is the \emph{support function} of $D$.}
It's easy to see that $\phi_D$ and $D$ determine each other. Let $\polytp_D = \{u \in M_{\R} ~ | ~ u \geq \phi_D ~\textrm{ on }~ N_{\R}\}$, then $\polytp_D$ is a (possibly empty) polytope and $\polytp_D \inter M$ generates the space of sections of the line bundle $\sheafo(D)$.
\prop{tDivisors}{\emph{(\cite{fulton} \S3.4)}
\begin{enumerate}
\item $D \stackrel{lin}{\sim} 0 \iff \phi_D$ is affine function,
\item $\linearsys{D}$ is generated by sections $\iff \phi_D$ is convex,
\item $D$ is ample $\iff \phi_D$ is strictly convex,
\item If $\polytp_D$ is non-singular, then $D$ is very ample $\iff \phi_D$ is strictly convex.
\end{enumerate}
In fact, all the relations among the divisors are generated by those of the following form
$$\sum_{\rho} \pairing{m}{v_{\rho}} D_{\rho} \stackrel{lin}{\sim} 0, ~\textrm{ sum over all }~ (n-1)\textrm{-faces of }~ \polytp, ~\textrm{ for some }~ m \in M.$$
}
We also note the following fact about canonical class of a smooth toric variety
\prop{canClass}{When $X_{\polytp}$ is smooth, the canonical class $K_X = -\sum_{\rho} D_{\rho}$, thus the support function $\phi_K$ of $K_X$ has value $-1$ on all the $v_{\rho}$'s, and $X_{\polytp}$ is Fano $\iff$ $\phi_K$ is strictly convex.}

When $\linearsys{D}$ is generated by sections, the polytope $\polytp_D$ also determines $\phi_D$, which in turn determines $D$. Thus $\polytp$ defines an ample divisor $D_{\polytp}$ on $X$. $D_{\polytp}$ is very ample if $\polytp$ is non-singular, and it defines an embedding of $X$ into $\CP^{\ell - 1}$, where $\ell = \numin{\polytp}$ is the number of lattice points in $\polytp$. This embedding is the same as the one given by $i_{\polytp}$ above.

\section{Semi-stable partition of polytopes}

This is a purely combinatoric section, in which we define (\emph{non-singular}, or \emph{mildly singular}) \emph{semi-stable} partition (\citedef{partition}, \citedef{balanced}) of a (non-singular, or simplicial) polytope $\polytp$, and construct a (non-singular, or simplicial) polytope $\new{\polytp}$ of one dimension higher from the partition. In next section, we'll show that $\new{\polytp}$ defines a toric variety, which is the total space of a degeneration of the toric variety defined by $\polytp$. The main results of this section are \citeprop{existence} and \citethm{construction}, and we won't need the details of the proofs in the rest of the paper.

Let $\Gamma$ be a partition of a simplicial polytope $\polytp$ into smaller simplicial polytopes $\polytp_j$'s. We'll call such partition \emph{simplicial} partition. An $l$-face $\face$ of $\polytp_j$ will also be called $l$-face of $\Gamma$ and be denoted by $\face \faceof \Gamma$. We make an exception for the $0$-faces of $\polytp$ by declaring they are \emph{not} a $0$-face of $\Gamma$. We'll call a $0$-face ($1$-face) a \emph{vertex} (\emph{edge}). The restriction of $\Gamma$ to $\face \faceof \polytp$ is defined to be the partition formed by $\{\polytp_j \inter \face\}$, and denote the restriction as $\Gamma \inter \face$, so the faces of $\Gamma \inter \face$ are also faces of $\Gamma$. We have the following definition for a simplicial partition:
\defin{partition}{$\Gamma$ is \emph{semi-stable}if for any $l$-face $\sigma$ of $\Gamma$ and $k$-face $\tau$ of $\polytp$, if $\sigma \subset \tau$, then there are exactly $k - l + 1$ $\polytp_j$'s such that $\sigma \faceof \polytp_j$.}
\exmp{partition1}{In dimension 1, where $M \equ \Z$ and $M_{\R} \equ \R$. $\polytp$ is the line segment between $0$ and $n \in \Z$. Any partition with partition points in $\Z$ is semi-stable.}
\exmp{partition2}{In dimension 2, where $M \equ \Z^2$ and $M_{\R} \equ \R^2$. Let $\polytp$ be the convex hull of the points $(0,0)$, $(3,0)$ and $(0,3)$. In the following partitions, the first two are semi-stable while the last one is not.
\begin{center}
\begin{tabular}{c c c}
\hspace{3.5cm} & \hspace{3.5cm} & \hspace{3.5cm} \\
\draw{
 \input{lattice.tex}
 \linewd 0.03 \move (0 0) \rlvec (0 1.5) \move (0 0) \rlvec (1.5 0) \rlvec (-1.5 1.5)
 \linewd 0.01 \move (0 1) \rlvec (1 -1)
} &
\draw{
 \input{lattice.tex}
 \linewd 0.03 \move (0 0) \rlvec (0 1.5) \move (0 0) \rlvec (1.5 0) \rlvec (-1.5 1.5)
 \linewd 0.01 \move (0.5 0.5) \rlvec (0 -0.5) \move (0.5 0.5) \rlvec (-0.5 0.5) \move (0.5 0.5) \rlvec (0.5 0)
} &
\draw{
 \input{lattice.tex}
 \linewd 0.03 \move (0 0) \rlvec (0 1.5) \move (0 0) \rlvec (1.5 0) \rlvec (-1.5 1.5)
 \linewd 0.01 \move (0.5 0) \rlvec (0 1) \move (0 1.5) \rlvec (0.5 -1)
} \\
(a) & (b) & (c) \\
\end{tabular}
\end{center}
}

The main goal of this section is to construct liftings of $\polytp$ by a semi-stable partition $\Gamma$.
\defin{lifting}{
A {\rm lifting} $(\new{\polytp}, \new{M}, \pi)$ of $\polytp$ by a semi-stable partition $\Gamma$ consists of an integral polytope $\new{\polytp}$ w.r.t. $\new{M}$ and a surjective $\pi : \new{M} \to M$, satisfying the following compatiblity condition. Let $\new{\face} \faceof \new{\polytp}$, $\pi_* : \new{M}_{\R} \to M_{\R}$ the map induced from $\pi$, then either $\pi_*(\new{\face}) \faceof \polytp$ or $\pi_*(\new{\face}) \faceof \Gamma$. If $\pi_*(\new{\face}) \faceof \Gamma$, $\new{\face}$ is said to be a \emph{lift} of $\pi_*(\new{\face})$ in $\new{\polytp}$. If all polytopes involved are non-singular, the lifting is called a \emph{non-singular lifting}.
}
We'll write the lifting as $(\new{\polytp}, \pi)$, or even $\new{\polytp}$ when the omitted is clear from the context. By the definition, the lifting will have dimension higher than that of $\polytp$. This section only deal with the liftings having dimension $\dim{\polytp} + 1$. We start by deriving some properties of semi-stable partitions. In the rest of the section, we'll assume all the partitions to be semi-stable.
\prop{inter}{
If $\inter_{k = 1}^{l} \polytp_{j_k} \neq \emptyset$, then it has dimension $n - l + 1$.

\proof{
First suppose $l = 2$. Assume the opposite, i.e., $\polytp_i \inter \polytp_j \neq \emptyset$ has dimension $l < n - 1$. Let $\tau$ be a $l$-face of $\Gamma$ that is contained in the intersection. Suppose $\tau$ lies in a $t$-face $\face_t$ of $\polytp$, then we have $t - l + 1$ $\polytp_k$'s that have $\tau$ as a face. Since $\polytp$ is simplicial polytope, there are $n - t$ $(n-1)$-faces of $\polytp$ which have $\face_t$ as a face. On the other hand, there are $n - l$ $(n-1)$-faces of $\polytp_i$ having $\tau$ as a face, and at most $n - t$ of them could lie on $(n-1)$-faces of $\polytp$, which left $t - l$ of them contained in the unique $n$-face of $\polytp$. By the assumption, none of the $(n-1)$-faces of $\polytp_i$ is also a $(n-1)$-face of $\polytp_j$. By the definition, there have to be $t - l$ $\polytp_k$'s different from $\polytp_i$ or $\polytp_j$ having $\tau$ as a face. In all, there are $t - l + 2$ $\polytp_k$'s having $\tau$ as a face, contradiction.

For general $l$, the $l - 1$ intersections $\polytp_{j_1} \inter \polytp_{j_k}$ $(k = 2, \ldots, l)$ are $(n-1)$-faces of $\polytp_{j_1}$, which has non-empty intersection. They have exactly one common $(n-l+1)$-face by property of simplicial polytope.
}
}
The lemma showed that two simplicial polytopes in a semi-stable partition can only intersect along the closure of an $(n-1)$-face.
Similarly, we have the following alternative lemma for intersection with faces of $\polytp$:
\lemma{interface}{
If $\polytp_i \inter \face_k \neq \emptyset$ for $\face_k \faceof \polytp$ a $k$-face, then it has dimension $k$.

\proof{
Suppose $\polytp_i \inter \face_k \neq \emptyset$. Let $\tau$ be a $l$-face of $\Gamma$ that is contained in the intersection. There are exactly $k - l + 1$ $\polytp_j$'s $\{\polytp_i, \polytp_{j_1}, \ldots \polytp_{j_{k-l}}\}$ that have $\tau$ as a face. By \citeprop{inter}, in the $n - l + 1$ $(n-1)$-faces of $\polytp_i$ that have $\tau$ as a face, $k-l$ of them are the intersections with $\polytp_{j_{t}}$'s, which left $n-k$ of those lie on the $n-k$ $(n-1)$-faces of $\polytp$ that have $\face_k$ as a face. So there is a $k$-face of $\polytp_i$ having $\tau$ as a face contained in $\face_k$.
}
}
From this lemma, we have that the restriction of a semi-stable partition to a face of $\polytp$ is a semi-stable partition of the face. As a corollary, we have:
\prop{interface}{
If $\inter_{i = 1}^{l} \polytp_{j_i} \inter \face_k \neq \emptyset$, for $\face_k \faceof \polytp$ a $k$-face, then it has dimension $k - l + 1$.

\proof{
Apply \citeprop{inter} to the restricted partition $\Gamma \inter \face_k$.
}
}
\prop{edgecount}{
$\Gamma$ is semi-stable then for any vertex $p$ of $\Gamma$ which is not a vertex of $\polytp$, there are exactly $n + 1$ edges $\face_0 \ldots, \face_n$ of $\Gamma$ such that $p \faceof \face_i$.

\proof{
First assume that $p$ is in the interior of $\polytp$, then from definition, there are $n - 0 + 1 = n+1$ $\polytp_j$'s having $p$ as a face. For each of the edges of $\Gamma$ that having $p$ as a face, there are $n - 1 + 1 = n$ of the $n+1$ $\polytp_j$'s intersect along it. It follows that the proposition is true in this case.

Suppose $p \subset \face$ for $\face$ a $k$-face of $\polytp$. By definition, there are $k - 0 + 1 = k + 1$ $\polytp_j$'s i.e., $\{\polytp_i, \polytp_{j_1}, \ldots \polytp_{j_k}\}$ which have $p$ as a face. Consider the partition $\Gamma \inter \face$, since $p$ is a interior point of $\face$, it's a face of $k + 1$ edges of $\Gamma \inter \face$. By \citeprop{inter}, the $\polytp_{j_t}$'s intersect along the closure of a $(n - k)$-face of $\Gamma$, which gives $n - k$ edges of $\Gamma$ that have $p$ as a face. Add the numbers, the proposition is proved.
}
}
From the above arguments, there is a simplicial complex $K_{\Gamma}$ associated to a semi-stable simplicial partition $\Gamma$. The vertices and simplices sets are defined as following:
\begin{itemize}
\item Vertices set is the set of polytopes in partition $\Gamma$,
\item For each $l$-face $\face$ in $\Gamma$ that lies in the interior of $\polytp$, an $(n - l)$-simplex $\simplex$ is defined by the polytopes that have $\face$ as a face. $\simplex$ and $\face$ are said to be the \emph{dual} of each other, and we also denote $\simplex$ as $\simplex_{\sigma}$.
\end{itemize}
\defin{dualcomplx}{The simplicial complex $K_{\Gamma}$ defined as above is the \emph{dual simplicial complex} of $\Gamma$.}
By definition, $K_{\Gamma}$ has a $k$-simplex if and only if there is a $k$-face of $\polytp$ which contains a vertex of $\Gamma$. It's also easy to see that $K_{\Gamma}$ is contractible since $\polytp$ is contractible.

Let $\Gamma$ be a semi-stable partition, and $p \faceof \Gamma$ be a vertex contained in an $l$-face $\tau$ of $\polytp$. Then there are $l + 1$ edges $\face_0, \ldots, \face_l$ of $\Gamma \inter \tau$ such that $p \faceof \face_i$ $(i = 0, \ldots, l)$. For $i = i, \ldots, l$, let $\face_i(p)$ denote the primary vector of $\face_i$ at $p$, then they satisfy a single integral relation $w_0 \face_0(p) + \ldots + w_l \face_l(p) = 0$. We may choose $w_i > 0$ for all $i$, and denote $\face'_i(p) = w_i\face_i(p)$. Reorder the $w_i$'s such that $w_1 \leq w_2 \leq \ldots \leq w_l$, and call the prime vector $W_p = (w_0, \ldots, w_l) \in \Z^{l+1}$ the \emph{weight vector} of $\Gamma$ at $p$.
\defin{balanced}{A vertex $p \faceof \Gamma$ is \emph{balanced} if $W_p = (1, \ldots, 1)$, $p$ is further \emph{non-singular} if $p$ is non-singular in one (thus all) sub-polytope containing $p$. A semi-stable partition is \emph{balanced} if all it's vetices are balanced, is \emph{non-singular} if furthermore all its vertices are non-singular. Suppose $K_{\Gamma}$ is $l$-dimensional, \emph{maximal vertices} of $\Gamma$ are the vertices of $\Gamma$ that lie in some $l$-face of $\polytp$. A balanced $\Gamma$ is \emph{mildly singular} if all its maximal vertices are non-singular.}

By definition, it's trivially true that if $p \faceof \Gamma$ is contained in an edge of $\polytp$, then $p$ is balanced. Also clear from the definition, any intersection among the subpolytopes of $\Gamma$ contains at least one maximal vertex. The following statement shows that the converse of the above definition for non-singular partition is also true, i.e., if all the vertices of a semi-stable partition $\Gamma$ are non-singular, then $\Gamma$ is balanced.
\prop{edgesum}{
Let $p \faceof \Gamma$ be a vertex contained in an $l$-face $\tau$ of $\polytp$. Let $\face_0, \ldots, \face_n$ be the $n+1$ edges of $\Gamma \inter \tau$ such that $p \faceof \face_i$ $(i = 0, \ldots, n)$, then $\sum_{i = 0}^l \face_i(p) = 0$.

\proof{
First suppose $l = n$, i.e. $\tau$ is the interior of $\polytp$. Since any $n$ of the $\face_i(p)$'s generate $M$ as a basis, we can write $\face_0(p) = \sum_{i = 1}^n v_i \face_i(p)$, $v_i \in \Z$ $(i = 1, \ldots, n)$. Let $e_i = (0, \ldots, 1, \ldots, 0) \in \Z^n$ with $1$ at the $i$-th place, and $v = (v_1, \ldots, v_n)$. The change of basis in $M$ between $\{\face_1(p), \ldots, \face_n(p)\}$ and $\{\face_0(p), \face_1(p), \ldots, \oomit{\face_i(p)}, \ldots, \face_n(p)\}$ is given by matrix $X_i = (\transpos{v} \transpos{e_1} \ldots \transpos{e_{i-1}} \transpos{e_{i+1}} \ldots \transpos{e_n})$. So $\det(X_i) = \pm 1$ for $i = 1, \ldots, n$, which gives $v = (\pm 1, \pm 1, \ldots, \pm 1)$. $p$ is a interior point of $\polytp$, so $v = (-1, \ldots, -1)$, which proves the statement.

For general $l$, the result follows from applying the above to $\Gamma \inter \tau$.
}
}

Let $\sigma \faceof \Gamma$ be a $(n-1)$-face, and suppose $\sigma = \polytp_i \inter \polytp_j$. Let $p \faceof \sigma$ be a vertex, then there are $n+1$ edges $\tau_0, \ldots, \tau_n$ of $\Gamma$ that meet at $p$. Suppose that $\tau_0 \nfaceof \polytp_i$, while $\tau_1 \nfaceof \polytp_j$. Let $f_{ij, p}$ be the rational affine function which defines the hyperplane that contains $\sigma$, and $f_{ij, p}(p + \tau'_0(p)) = 1$. By definition, $f_{ij, p} = -f_{ji, p}$. For another vertex $q \faceof \sigma$, we can define $f_{ij, q}$ similarly, and it will differ from $f_{ij, p}$ in general. When $p$ and $q$ are both non-singular, however, we have $f_{ij, p} = f_{ij, q}$.

Suppose $K_{\Gamma}$ is $l$-dimensional. When $l = 1$, $\Gamma$ is balanced, and $\Gamma$ partitions $\polytp$ into a chain of subpolytopes by a set of non-intersecting hyperplanes. Choose a vertex $p_{ij}$ for the partition hyperplane between $\polytp_i$ and $\polytp_j$, and define $f_{ij} = f_{ij, p_{ij}}$. For $l > 1$, we'll assume that $\Gamma$ is mildly singular for the rest of the section. We can define $f_{ij}$ using only the maximal vertices, which are all non-singular under the assumption we made.

Let $\family{A} = \{f : \R^n \to \R ~|~ f \textrm{ affine, and } f(\Z^n) \subset \Q\}$ be the group of all rational affine functions on $\R^n$. $\family{A}$ can be identified with $\Q^{n+1}$. Let $C(K_{\Gamma})$ be the oriented chain complex of $K_{\Gamma}$ (\cite{spanier} p.170), then we define $\alpha \in Hom(C_1(K_{\Gamma}), \family{A})$ by $\alpha(\simplex_{\polytp_i \inter \polytp_j}) = f_{ij}$ on the basis and extend linearly to the whole $C_1(K_{\Gamma})$.
\lemma{cocycle}{ The cochain $\alpha$ as defined above is a cocycle.

\proof{
Suppose $K_{\Gamma}$ is $l$-dimensional. For $l = 1$ there is nothing to prove. We'll assume $l >1$.

Without lost of generality, we assume that $i$, $j$, $k$ are mutually distinct, $p \faceof \polytp_k$ is non-singular, and $\tau_2 \nfaceof \polytp_k$. Since $p$ is non-singular, $\tau'_i(p) = \tau_i(p)$ for all $i$. Let $\face_{n-2} = \polytp_i \inter \polytp_j \inter \polytp_k$, then we have the following:
$$f_{ij} = f_{jk} = f_{ki} = 0 \hspace{1in} {\rm on} \,\, \face_{n-2}$$
$$\equlist{lll}{
f_{ij}(p + \tau_0(p)) = 1, & f_{ij}(p + \tau_1(p)) = -1, & f_{ij}(p + \tau_2(p)) = 0 \\
f_{jk}(p + \tau_0(p)) = 0, & f_{jk}(p + \tau_1(p)) = 1, & f_{jk}(p + \tau_2(p)) = -1 \\
f_{ki}(p + \tau_0(p)) = -1, & f_{ki}(p + \tau_1(p)) = 0, & f_{ki}(p + \tau_2(p)) = 1
}$$
which shows $f_{ij} + f_{jk} + f_{ki} = 0$.
}
}
Since $K_{\Gamma}$ is contractible, the cohomology group $\coho{1}{K_{\Gamma}, \family{A}} = 0$. It follows that $\exists \beta \in Hom(C_0(K_{\Gamma}), \family{A})$, s.t. $\cod(\beta) = \alpha$. Specifically, $\beta$ assigns a function $f_i \in \family{A}$ to each sub-polytope $\polytp_i$ in $\Gamma$, such that $f_j - f_i = f_{ij}$ if $\polytp_i \inter \polytp_j \neq \emptyset$. Suppose $\cod{\beta'} = \alpha$ as well, let $\gamma = \beta - \beta'$ then $\cod{\gamma} = 0$, i.e. $\gamma$ assigns the same function $f \in \family{A}$ to all the sub-polytopes in $\Gamma$. We can define a function $F$ on $\polytp$ by $F \restricto{\polytp_i} = f_i$. We'll use the notion \emph{piecewise affine function on $\Gamma$} to denote a piecewise affine function on $\polytp$ s.t. it's affine on each face of $\Gamma$. Then $F$ is a piecewise affine function on $\Gamma$.

\defin{concavity}{Let $f$ be a piecewise affine function on $\Gamma$, $p \faceof \Gamma$ be a vertex in the interior of $\polytp$ and the $n+1$ edges that meet at $p$ be $\face_0, \ldots, \face_n$.  Let $\delta_i(f, p) = f(p + \face_i(p)) - f(p)$ for $i = 0, \ldots, n$. $\concav(f, p) = \sum_{i = 0}^n \delta_i(f, p)$ is the \emph{concavity} of $f$ at $p$. If $p \in \face \faceof \polytp$ where $\face$ is a $0 < k$-face, we define $\concav(f, p) = \concav(f\restricto{\face}, p)$.}
Note the concavity is not defined at the vertices of $\polytp$. $\concav(-, p)$ is linear, i.e. $C(af + bg, p) =a C(f, p) +b C(g, p)$ for $a , b \in \R$. $f$ is a global affine function on $\Gamma$, iff $C(f, p) = 0$ for all $p$ where it is defined. It implies that if $C(f, p) = C(g, p)$ for all $p$, then $f - g$ is a global affine function. $f$ is convex (concave) at $p$ (i.e. in a neighbourhood of $p$) iff $C(f, p) < 0$ (resp. $> 0$). So as the name suggested, $C(f, p)$ captures the concavity of $f$ at $p$. We see that $F$ as constructed above has concavity $0 < C(F, p) \in \Q$ at all $p$, so is concave. Summing up, we have

\prop{existence}{
For a mildly singular semi-stable partition $\Gamma$, we can define a concave rational piecewise affine function $F_{\Gamma}$ on $\Gamma$. We call such $F_{\Gamma}$ a \emph{(rational) lifting function} of $\Gamma$. Lifting function of $\Gamma$ is unique modulo global affine functions and multiplication by rationals.
\eopf
}
Let $C(F_{\Gamma}) = \{C(F_{\Gamma}, p) ~|~ 0\textrm{-face } p \faceof \Gamma\}$, then in general, $C(F_{\Gamma})\subset \Q$. For balanced $\Gamma$, choices can be made such that $C(F_{\Gamma}) = \{ 1 \}$. Multiply $F_{\Gamma}$ by $r \in \Q$, we get $F'_{\Gamma} = r F_{\Gamma}$, with $C(F'_{\Gamma}) = r C(F_{\Gamma})$. For suitable $r$, $F'_{\Gamma}$ is integral in the sense $F'_{\Gamma}(\Z^n) \subset \Z$. Let $R$ be the minimal $r$ such that $F'_{\Gamma}$ is integral, and call $F'_{\Gamma}$ the \emph{minimal integral lifting}.

Suppose $F$ is a minimal integral lifting function of $\Gamma$. Let $\new{\polytp} = \{(y, x) \, | \, y \geq F(x)\} \subset \R \cross \polytp \subset \R \cross M_{\R}$.
\thm{construction}{
Let $\pi$ be the projection $\Z \dsum M \to M$, then $(\new{\polytp}, \Z \dsum M, \pi)$ is a lifting of $\polytp$ by $\Gamma$, and there is exactly one lift of each $\face \faceof \Gamma$ in $\new{\polytp}$. If $\Gamma$ is a non-singular partition of a non-singular polytope, and $C(F) = \{1\}$, then the lifting is non-singular.

\proof{
Everything basically follows from the construction.

$\new{\polytp}$ is convex and integral because $F$ is concave and piecewise integral affine. It's simplicial by \citeprop{edgecount} and the fact that there is a vertical edge over each vertex of $\polytp$. $\new{\face} \faceof \new{\polytp}$ is either  $\{(y, x) \in \R \cross \face \, | \, y \geq F(x)\}$ for the vertical faces or the graph of a integral affine function on $\face = \pi_*(\new{\face})$. It follows that either $\face \faceof \polytp$ or $\face \faceof \Gamma$, and each face of $\Gamma$ has exactly one lift in $\new{\polytp}$. Nonsingularity follows from the concavity condition.
}
}

Sometimes it's more convenient to work with compact polytopes instead of open ones. In that case, we can consider $\new{\polytp}^c = \new{\polytp} \inter \{(y, x) \, | \, y \leq \pairing{a}{x} + b\}$ for $a \in N$, and sufficiently big $b \in \Z$. We call $\new{\polytp}^c$ a \emph{compact} lifting of $\polytp$ by $\Gamma$. It can be thought of as adding a face to $\new{\polytp}$ at infinity.

\section{Semi-stable degeneration of toric varieties}

We'll prove the main result (\citethm{degeneration}) and give some examples. There is also a similar result (\citethm{weakdegen}) on the degeneration of toric varieties with orbifold singularities. As a result of the construction, we give concrete description of the degenerations in \citermk{degen3}. When we talk about polytopes and partitions, we'll use the notations introduced in last section.

Let $\polytp \subset M_{\R}$ be a simplicial polytope, $\fan$ be its fan, and the corresponding toric variety be $X$. Suppose $\Gamma$ is a non-singular semi-stable partition of $\polytp$. Recall that the dual graph of the central fiber of a semi-stable degeneration is defined by assigning a vertex for each component, and a $k$-simplex whenever $k$ components having non-empty intersection.
\thm{degeneration}{Suppose $\polytp$ is non-singular, then there exists a semi-stable degeneraiton $p : \new{X} \to \C$ of $X$ to $\inverse{p}(0)$. The dual graph $G$ of the fiber $\inverse{p}(0)$ is isomorphic to $K_{\Gamma}$, and a component in $\inverse{p}(0)$ is the toric variety defined by the corresponding subpolytope in $\Gamma$.

\proof{
By \citethm{construction}, let $\new{\polytp}$ be a non-singular lifting of $\polytp$ by $\Gamma$. Let $\new{\fan}$ be the fan of $\new{\polytp}$, then $\fan_{\Gamma}$ is a subfan of $\new{\fan}$. Let $\new{X}$ be the toric variety defined by $\new{\fan}$, then $\new{X}$ is smooth.

Note that by construction, the support of $\new{\fan}$ is the closed upper halfspace $\cl{\upperhalf} = N_{\R} \cross \cl{\Rplus}$. The cones $\cone_{\face}$ that correspond to the vertical faces $\face$ of $\new{\polytp}$ lie in $N_{\R} \cross \{0\}$, and all the other cones in $\new{\fan}$ lie in $N_{\R} \cross \Rplus$. Let $\Lambda$ be the fan in $\R$ with $\support{\Lambda} = \cl{\Rplus}$, then the variety it defines is $\C$.

Consider the exact sequence of lattices:
\equlabel{lattseq}{\shortexact{0}{N}{i}{N \dsum \Z}{\mu}{\Z}{0},}
where $\mu$ is the projection to the second summand and $i$ is the inclusion. The maps are compatible with the fans $\fan$, $\new{\fan}$ and $\Lambda$, and define a sequence of proper morphisms:
\equlabel{morphseq}{\shortexact{pt}{X}{i_*}{\new{X}}{p}{\C}{pt},}
from which it follows that the fiber of $p$ over the generic point of $\C$ is isomorphic to $X$. The rays of $\new{\fan}$ in $\upperhalf$ is mapped to the open positive cone of $\R$, thus $p$ maps the divisors defined by such rays to $0$, i.e. $\inverse{p}(0) = \union_{i = 1}^k X_i$, where $X_i$'s are the divisors defined by such rays. By construction, $X_i$'s are themselves toric varieties which is defined by the subpolytopes in the partition $\Gamma$. The statement about the dual graph of the fiber $\inverse{p}(0)$ is obvious.

Next we write down the map $p$ in the affine pieces. Dualizing \citeequ{lattseq}, we have exact sequence: \equlabel{dualseq}{\shortexact{0}{\Z}{j}{M \dsum \Z}{\nu}{M}{0},} where $j$ is the inclusion, and $\nu$ is the projection to the first summand. Suppose $\new{\polytp}$ is defined by the lifting function $F$. Let $q \faceof \Gamma$ be a vertex which is an inner point of $\polytp$, and $\new{q}$ be the lift of it in $\new{\polytp}$. Let $\face_0, \ldots, \face_n$ be the $n + 1$ edges of $\new{\polytp}$ that meet at $\new{q}$, and $\cone_{q}$ be the cone dual to $\new{q}$, then $\dualcn{\cone_{q}}$ is generated by $\face_0(\new{q}), \ldots, \face_n(\new{q})$. Let $X_{q}$ be the affine piece defined by $\dualcn{\cone_{q}}$. Since $\Gamma$ is non-singular, we have $C(F) = \{1\}$. It follows that $\sum_{k = 0}^{n}\face_k(\new{q}) = (0, \ldots, 0, 1) = j(1)$, so the map $j$ induces $j_* : \C[t] \to \C[x_0, \ldots, x_n], t \mapsto \prod_{k = 0}^{n} x_k$, where $x_k$'s are the coordinates given by the edges of $\dualcn{\cone_{q}}$. It follows that the map $p$ on the affine piece $X_{q}$ is defined by $(x_0, \ldots, x_n) \mapsto \prod_{k = 0}^{n}x_k$.

Now suppose $q \faceof \Gamma$ is a vetex which is an inner point of a $0 < r$-dimensional face $\tau \faceof \polytp$. By \citeprop{edgesum}, there are $r+1$ edges of $\Gamma \inter \tau$ meeting $q$, and they add to $0$. Let $\face_0, \ldots, \face_r$ be the corresponding $r+1$ edges of $\new{\polytp}$ that meet $\new{q}$. Then $\sum_{k = 0}^{r}\face_k(\new{q}) = (0, \ldots, 0, 1) = j(1)$ by concavity, and $j_*$ is defined by $t \mapsto \prod_{k = 0}^r x_k$, and $p$ is defined by $(x_0, \ldots, x_n) \mapsto \prod_{k = 0}^{r}x_k$ in the affine piece $X_q$.

If $q$ is a vertex of $\polytp$, then it is a face of a vertical edge $\face_0$ of $\new{\polytp}$. Thus $j_*(t) = x_0$, and $p(x_0, \ldots, x_n) = x_0$ in the affine piece $X_q$.

Since the affine pieces as described above cover the central fiber, the theorem follows.
}
}
\rmk{degen2}{If we take a compact lifting, then instead of \citeequ{morphseq}, we get sequence of morphisms of compact toric varieties:
\equlabel{morphseq2}{\shortexact{pt}{X}{i_*^c}{\new{X}^c}{p^c}{\CP^1}{pt},}
}

\rmk{degen3}{We can write down the degeneration in more concrete terms using the embedding of $X$ into $\CP^{\ell-1}$ given by $\polytp$ (recall that $\ell = \numin{\polytp}$). Suppose $F$ and $F'$ are lifting functions that give rise to non-singular liftings of $\polytp$ as in \citethm{construction}, then $F' = F+G$ for some integral affine function $G$.
Let $\new{i}_{F} : (\ctori)^{n+1} \to \CP^{\ell - 1} \cross \C$ be the embedding defined by
$$\new{x} = (x, \lambda) \mapsto ([\lambda^{F(\lat{m}_1)}x^{\lat{m}_1} : \ldots : \lambda^{F(\lat{m}_{\ell})}x^{\lat{m}_{\ell}}], \lambda).$$
Let $\new{i}_{F'}$ be the similar embedding defined by $F'$. Then $\new{i}_{F}$ and $\new{i}_{F'}$ has the same images since $\lambda \neq 0$. Let $\new{X}_F$ be the closure of the image of $\new{i}_F$.  It follows that $\new{X} = \new{X}_{F}$ is independent of the lifting, and it's the same $\new{X}$ as in \citethm{degeneration}. The map $p$ is given by
$$p : ([\lambda^{F(\lat{m}_1)}x^{\lat{m}_1} : \ldots : \lambda^{F(\lat{m}_{\ell})}x^{\lat{m}_{\ell}}], \lambda) \mapsto \lambda.$$
Suppose $\Gamma$ partition $\polytp$ into $\union_{i=0}^k \polytp_i$, and let $X_i$ be the toric variety defined by $\polytp_i$. By adding an affine function, we may assume that $F\restricto{\polytp_0} = 0$, then $F(\lat{m}) > 0$ for $\lat{m} \notin \polytp_0$ by concavity.
Let $[y_1 : \ldots : y_{\ell}]$ be the homogeneous coordinates of $\CP^{\ell -1}$, and define
$$H_i = \{[y_1 : \ldots : y_{\ell}] | y_{j} = 0 \textrm{ if } \lat{m}_j \notin \polytp_i\} ~~~~\textrm{ for } ~ i = 0, \ldots, k.$$
Let $\ell_i = \numin{\polytp_i}$, then $H_i \iso \CP^{\ell_i - 1}$. Taking limit of the embedding $\new{i}_{F}$ as $\lambda \to 0$, we get the embedding of $X_0$ into $\CP^{\ell_0 - 1} \iso H_0$. Using different presentations of $\new{X}$ by using different lifting functions, and taking limit as $\lambda \to 0$, we see that $\inverse{p}(0) = \union_{i = 0}^k X_i$.
}

In case of $\polytp$ being singular, we can still formulate a similar result using the concept of ``weak semi-stable degeneration'' as in \cite{ambramovich}. The difference between weak semi-stable and semi-stable degenerations is that the total space in the first case need not to be non-singular. The degeneration we construct from a singular simplicial toric variety could not have non-singular total space, since $\polytp$ is not non-singular, while we can still require the partition to be mildly singular.
\thm{weakdegen}{Suppose $\Gamma$ is a mildly singular semi-stable partition of simplicial $\polytp$, then there exists a weak semi-stable degeneration $p : \new{X} \to \C$ of $X$ to $\inverse{p}(0)$. The dual graph $G$ of the fiber $\inverse{p}(0)$ is isomorphic to $K_{\Gamma}$, and a component in $\inverse{p}(0)$ is the toric variety defined by the corresponding subpolytope in $\Gamma$.

\proof{
The proof goes word by word as the proof of \citethm{degeneration}, except for the fact that $\new{X}$ might be singular. The reducibility of the central fiber follows from the fact that $\Gamma$ is balanced.
}
}

\exmp{degenCP1}{The semi-stable degeneration for the partition in \citexmp{partition1} is the result of successively blowup smooth points in the central fiber of $\CP^1 \cross \C$.}

\exmp{degenCP2}{In \citexmp{partition2}, $(a)$ corresponds to blowup one point in the central fiber of $\CP^2 \cross \C$. The proper transformation of the central fiber is then $\blowup{\CP^2} \iso \p(\bundle{O}(-1) \dsum \bf{1})$. Blowing-up a fiber in $\p(\bundle{O}(-1) \dsum \bf{1})$, we get the degeneration corresponding to $(b)$.}

\exmp{degenCPn}{The construction as in \citexmp{degenCP2} can be generalized to $\CP^n$. Let $e_0 = (-1, -1, \ldots, -1) \in \Z^n \subset \R^n$, and $\{e_i\}$ be the standard integral basis of $\Z^n$. Let $\polytp(n)$ be the polytope in $\R^n$ spanned by $e_0$ and $e_0 + (n+1)e_i$, $i = 1, \ldots, n$, then $\polytp(n)$ is the standard reflexive polytope that represents $\CP^n$. Let $v_0 = e_1$, $v_i = e_{i+1} - e_{i}$ for $i = 1, \ldots, n$, where $e_{n+1} = 0$. Let $\fan(n)$ be the maximal fan generated by $\{v_i\}_{i = 0}^{n}$, then $\fan(n)$ is complete and nonsingular. $\polytp(n) \inter \fan(n)$ gives a non-singular semi-stable partition $\Gamma(n)$ of $\polytp(n)$. There are $(n+1)$ subpolytopes in $\Gamma(n)$, and let $\polytp_i(n)$ denote the subpolytope containing $e_i$. The equations of the hyperplanes that bound $\polytp_0(n)$ are:
\equlabel{bounding}{\sum_{i = 1}^{j} x_i = 0, ~ \textrm{ for } ~ j =
1, \ldots, n.}
The following is the partition $\Gamma(3)$ for $\CP^3$:

\begin{center}
\draw{
 \linewd 0.05
 \rlvec (4 0) \rlvec (-4 4) \rlvec (0 -4)
 \linewd 0.01
 \rlvec (1.732 1) \rlvec (2.268 -1) \move (1.732 1) \rlvec (-1.732 3)
 \move (1 2) \linewd 0.03
 \lvec (0 3) \linewd 0.02 \lvec (0.866 1.5) \lvec (1.299 1.75)
 \lvec (2.433 1.25) \linewd 0.03 \lvec (2 2) \lvec (1 2)
 \lvec (1.433 1.25) \lvec (1.433 0.25)
 \linewd 0.02
 \lvec (0.866 0.5) \lvec (0.866 1.5) \lvec (1.433 1.25)
 \linewd 0.03 \lvec (2.433 1.25) \lvec (3.433 0.25)
 \lvec (1.433 0.25) \lvec (1 0)
 \lvec (1 2)
}
\end{center}

The semi-stable degeneration constructed from $\Gamma(n)$ will then have $n+1$ components in the central fiber, say $C_0, \ldots, C_n$. The degeneration can also be viewed as the result of blow-up smooth subvarieties in the central fiber of the trivial family as following. Let $\flag$ be a complete flag in $\CP^n$, i.e. a chain of linear subspaces $pt = \flag_0 \subsp \flag_1 \subspace \ldots \subsp \flag_{n-1}$ where $\flag_i$ has dimension $i$. Blow-up $\flag_i$ in the central fiber starting from $i = 0$, and then blow-up the proper transformation of $\flag_{i+1}$, until all the $\flag_i$'s are blown-up. By symmetry of the partition, we see that all the components in the resulting central fiber are isomorphic to $Q^n$, the proper transformation of $\CP^n$ after all the blow-ups. Below is an illustration of the polytopes for $Q^i$, $i = 1, 2, 3$:

\begin{center}
\begin{tabular}{c c c}
\hspace{3.5cm} & \hspace{3.5cm} & \hspace{4.5cm} \\
\draw{
 \linewd 0.03
 \rlvec (1 0)
} &
\draw{
 \linewd 0.03
 \rlvec (2 0) \rlvec (-1 1) \rlvec (-1 0) \rlvec (0 -1)
} &
\draw{
 \linewd 0.03
 \rlvec (2 0) \rlvec (0.354 0.354)
 \linewd 0.01
 \rlvec (-1.293 0.707) \rlvec (-1.061 -1.061)   
 \move (1.061 1.061) \rlvec (-0.354 0.646) \move (0 0)
 \linewd 0.03
 \rlvec (0 1) \rlvec (1 0) \rlvec (0.354 0.354) \rlvec (-0.647 0.353) \rlvec (-0.707 -0.707)   
 \move (2 0) \rlvec (-1 1) \move (2.354 0.354) \rlvec (-1 1)
} \\
$Q^1 \iso \CP^1$ & $Q^2 \iso \blowup{\CP^2}$ & $Q^3$
\end{tabular}
\end{center}
}

Of course, we may choose to partition $\polytp(n)$ into fewer pieces, which corresponding to either use some subset of the $\{v_i\}$ to generate a fan that has fewer cones, or blowup an incomplete flag in $\CP^n$. The resulting degeneration will correspondingly have fewer components in the central fiber.

\exmp{degenWPn}{We can similarly define a weak semi-stable degeneration of weighted projective space. Let $\polytp$ be the polytope in $\R^4$ spanned by $e_0$, $e_0 + 8 e_1$ and $e_0 + 4 e_i$ for $i = 2, 3, 4$, then $\polytp$ is the reflexive polytope representing the weighted projective space $\p(1, 1, 2, 2, 2)$. We use the fan $\fan(4)$ as described above to give a semi-stable partition of $\polytp$. The resulting partition is no longer symmetric, since $\polytp$ is not symmetric to begin with. It's not hard to show that the partition is balanced. It's mildly singular because the only maximal vertex $(0, 0, 0, 0)$ is non-singular.
}

\exmp{gendegen1}{It is not true that all the semi-stable degeneration constructed from the theorem can be constructed also by blowing up the central fiber of a trivial family. The following partition is a semi-stable partition, while the degeneration is not composition of blowing-ups on the central fiber.
\begin{center}
\draw{
 \input{latticebg.tex}
 \linewd 0.03 \move (0 1) \rlvec (0.5 -0.5) \rlvec (1 0) \rlvec (0.5 0.5)
 \rlvec (0 0.5) \rlvec (-0.5 0.5) \rlvec (-1 0) \rlvec (-0.5 -0.5) \rlvec (0 -0.5)
 \linewd 0.01 \move (1 0.5) \rlvec (0 1.5)
}
\end{center}
}

\section{Semi-stable degeneration of toric hypersurfaces}

Let $p : X \to \C$ be a semi-stable degeneration, and $V \subset X$ be a subvariety of $X$. Suppose $V$ intersects transversely with the central fiber. Let $p_V = p\restricto{V} : V \to \C$ and $V_{\lambda} = \inverse{p_V}(\lambda)$. Then over a small neighbourhood $U \subset \C$ of $0$, $p_V$ is a semi-stable degeneration of $V_{\lambda}$. The construction in the last section provides a nice way to get semi-stable degeneration of toric hypersurfaces and complete intersections. 

The subsections in this sections are organized as following. Subsection $5.1$ describes the two basic constructions to obtain semi-stable degenerations of toric subvarieties from the degeneration constructed in previous section. The construction with open polytopes is more explicit, while the construction with compact polytopes is a lot more general. Subsection $5.2$ proves a lemma which is quoted in subsequent subsections to determine the components in the singular fiber of the degenerations. The examples start with subsection $5.3$, in which we construct the semi-stable degeneration of degree $d$ surfaces in $\CP^3$ into a chain of rational components. In the same subsection, we obtain an algebraic geometrical version of the fiber-sum operation of two copies of $E(1)$ to get $K3$ surface. In subsection $5.4$, we recover the classic semi-stable degeneration of $K3$ surfaces into four rational components, and generalize it to degenerations of higher dimensional Calabi-Yau hypersufaces in projective spaces (\citethm{CYdegen}). Also in the same subsection, we briefly describe a weak semi-stable degeneration of Calabi-Yau hypersurface of $\p(1, 1, 2, 2, 2)$. The degenerations described in $5.3$ can be easily generalized to higher dimensions, and it's described in subsection $5.5$.

\subsection{Basic constructions}

\subsubsection{Construction with open polytopes} 
Let $\polytp \subset \R^n$ be a non-singular polytope, and $X$ be the toric variety defined by $\polytp$. $\polytp$ also defines a very ample divisor $D_{\polytp}$ in $X$. Let $\Gamma$ be a non-singular semi-stable partition of $\polytp$. We first consider a generic hypersurface $V$ in the linear system $\linearsys{D_{\polytp}}$. Referring to \citermk{degen3}, $V$ is given by the intersection of the image of $i_{\polytp}$ with a generic hyperplane in $\CP^{\ell-1}$, where $\ell = \numin{\polytp}$ is the number of lattice points in $\polytp$. We can thus write down the equations of $V_{\lambda} \subset X$ over a nonzero $\lambda$. Fix a generic hyperplane in $\CP^{\ell - 1}$:
$$H = \{[y_1: \ldots: y_{\ell}] | \sum_{j = 1}^{\ell} a_j y_j = 0\}.$$
The intersection of the image of $\new{i}$ with $H$ as subvarieties in $\CP^{\ell - 1} \cross \C$ gives the total space $\new{V}$ of the degeneration. Thus the equation of $V_{\lambda}$ is given by:
$$\phi_{\lambda} = \sum_{j = 1}^{\ell} \lambda^{F(\lat{m}_j)} a_j x^{\lat{m}_j} = 0.$$
We may again assume that $F\restricto{\polytp_0} = 0$ for some subpolytope $\polytp_0$. Let $X_0$ (and $D_0$) be the toric variety (resp. the divisor on $X_0$) defined by $\polytp_0$. Taking limit when $\lambda \to 0$, we see that the central fiber $V_0$ contains as a component a generic hypersurface of $X_0$ in the linear system $\linearsys{D_0}$.

\subsubsection{Construction with compact polytopes}
We may consider more general hypersurfaces and complete intersections by looking at compact lifting $\new{X}^c$ of $\polytp$ by $\Gamma$. Then by \citeequ{morphseq2}, $p^c$ gives a family over $\CP^1$. We'll drop the superscript ${}^c$ in the following. Suppose $\new{D}$ is a divisor on $\new{X}$, s.t. $\linearsys{\new{D}}$ is generated by sections. Then, a generic hypersurface in $\linearsys{\new{D}}$ is smooth, and it intersects transversely with generic fiber of $p$. Over a neighbourhood of $0 \in \CP^1$, we have a semi-stable degeneration of $V'_{\lambda} = V' \inter \inverse{p}(\lambda)$. We use $\rho$ to denote the vertical faces in $\new{\polytp}$, $\sigma$ to denote the faces from the partition $\Gamma$, and $\tau$ to denote the face at infinity, and we'll abuse the notation to denote the faces in $\polytp$ as $\rho$ as well. Then $D_{\tau}$ is the class of a generic fiber. Suppose $\new{D} = \sum_{\rho} a_{\rho} D_{\rho} + \sum_{\sigma} a_{\sigma} D_{\sigma} + a_{\tau} D_{\tau}$ and let $D = \sum_{\rho} a_{\rho} D_{\rho}$, then $V'_{\lambda}$ is in the linear system $\linearsys{D}$ of $X$. There is an obvious relation in $A_{n}(\new{X})$:
$$D_{\tau} = \sum_{\sigma} D_{\sigma}.$$
Suppose $\{\sum_{\rho} \alpha_{\rho}(i) D_{\rho} = 0 | i = 1, \ldots, n\}$ are the relations in $A_{n-1}(X)$. Then, besides the obvious one, all the relations in $A_{n}(\new{X})$ are given in 
$$\{\sum_{\rho} \alpha_{\rho}(i) D_{\rho} = \sum_{\sigma} \beta_{\sigma}(i) D_{\sigma} | i = 1, \ldots, n \} ~\textrm{ for some }~ \beta_{\sigma}(i).$$
Using these relations, we can get information on the components $C_i$ in the central fiber $V'_0$. It's easy to see that we can talk about complete intersections in this way as well.

We can write down the degeneration in more concrete terms, using the homogeneous coordinates for toric varieties (cf. \cite{cox}). Let $x_{\sigma}$ (resp. $x_{\rho}$, $x_{\tau}$) be the variables corresponding to the faces $\sigma$ (resp. $x_{\rho}$, $x_{\tau}$), then they have degrees $D_{\sigma}$ (resp. $D_{\rho}$, $D_{\tau}$) in the homogeneous ring $S(\new{X})$. The morphism $p : \new{X} \to \CP^1$ can be given by a morphism between the homogeneous rings $\C[z_0, z_1] \to S(\new{X})$, and $S(\new{X})$ is a graded $\C[z_0, z_1]$-algebra under this map. A complete intersection $V' \subset \new{X}$ is defined by the ideal $\ideal{I}$ generated by $\{f'_1, \ldots, f'_k ~ | ~ f'_i \in S(\new{X})$ is homogeneous of degree $ D'_i \in A_n(\new{X}) $ for $ i = 1, \ldots, k \}$. Given $\ideal{I}$, it's easy to find the ideal in $S(X)$ that defines a generic fiber $V'_{\lambda}$. To find the component in the central fiber that lies in $D_{\sigma}$, we look at the ideal in $S(\new{X})$ generated by $\ideal{I}$ and $x_{\sigma}$.

\subsection{A Lemma}

Before the examples, we need a lemma to identify the components in the central fibers. Suppose $X$ is a smooth projective variety of dimension $n$. Let $L$ be a Cartier divisor on $X$, and $\bundle{L}$ be the corresponding line bundle over $X$, with a generic section $s$. Let $\pi: E = \p_{X}(\bundle{L} \dsum \triv) \to X$ be the projectivization of $\bundle{L}$, and $f$ be the class of the fiber. Let $\lambda = [\lambda_0: \lambda_1] \in \CP^1$ then $(\lambda_1 s, \lambda_0)$ are sections of $\bundle{L} \dsum \triv$. Let $V_0$ (resp. $V_{\infty}$, $V_{\lambda}$) be (the closure of) the image of $(0, 1)$ (resp. $(s, 0)$, $(\lambda_1 s, \lambda_0)$) in $E$. $V_{\lambda}$ are sections of the $\CP^1$ bundle $E$. As divisors in $E$:
$$V_{\lambda \neq [1:0]} = V_{\infty} + \inverse{\pi}(L), ~~~\textrm{ and }~ V_0 \restricto{V_0} = L.$$
Let $D$ be an ample divisor in $X$ such that $L + D$ is also ample, and $\bundle{D}$ be the corresponding line bundle. Let $s_D$ (resp. $s'_D$) be a generic section of $\bundle{D}$ (resp. $\bundle{L} \tensor \bundle{D}$), and $\bundle{L}_D = (\bundle{L} \dsum \triv) \tensor \bundle{D} = \bundle{L} \tensor \bundle{D} \dsum \bundle{D}$, then $s_D$ (resp. $s'_D$) defines a generic hypersurface $D_{\infty}$ (resp. $D_0$) in the linear system $\linearsys{D}$ (resp. $\linearsys{L+D}$). Let $T = D_0 \inter D_{\infty}$, then $T$ is smooth codimension $2$ subvariety of $X$ since everything is generic. Consider the sections $s_{\lambda, D} = (\lambda_1 s'_D, \lambda_0s_D)$ of $\bundle{L}_D$. Let $W_{\lambda}$ be the closure of the image of $s_{\lambda, D}$ in $E \iso \p_X(\bundle{L}_D)$, then for $\lambda \neq 0, \infty$,
\equlabel{inter}{W_{\lambda} \inter V_0 = D_0, ~~ W_{\lambda} \inter V_{\infty} = D_{\infty}, ~\textrm{ and }~ W_{\lambda} \cdot f = 1.}
Apply Theorem IV-23 and Proposition IV-25 in \cite{eh00}, we see that $\pi\restricto{W_{\lambda}} : W_{\lambda} \to X$ is the blow-up of $X$ along $T$ for $\lambda \neq 0, \infty$. Denote the blow-up as $\blowup{X}_T$.

We claim that \citeequ{inter} essentially characterizes the blow-up of $X$ along $T$ as divisor in $E$.
\lemma{blowup}{Suppose $R$ is a divisor in $E$ that satisfies \citeequ{inter}, where the first two equations will be replaced by the corresponding intersections of divisor classes. Then a generic hypersurface $W$ in the linear system $\linearsys{R}$ is the blow-up of $X$ along $T$, where $T$ now is interpreted as intersection of generic hypersurfaces of $X$ in the linear systems $\linearsys{D}$ and $\linearsys{D+L}$.

\proof{ By \citeequ{inter}, $W$ is generically a section of $E$ over $X$, and $W$ has at least $2$ points over the points in $T$. Thus $W$ contains the whole fibers over points in $T$. The restriction of $\pi$ to $W$ is surjective and $\inverse{\pi}(T) = E\restricto{T}$ is Cartier in $W$. By the universal property of blow-up, $\pi\restricto{W}$ factors through a unique morphism $\phi: W \to \blowup{X}_T$. Well, $\phi$ has to be isomorphism since $W$ and $\blowup{X}_T$ have identical exceptional loci.
}
}
For the special case of $L = 0$, $E = X \cross \CP^1$, this gives the construction of a Lefschetz pencil. The $W$ in the lemma is the closure of the graph of a map from $X \sminus T$ to $\CP^1$.
\rmk{blowup1}{The $T$ in the lemma is \emph{not generic}, since it has to be the intersection of two hypersurfaces. For example, let $X = \CP^2$, $D = 2H$, $L = H$, then $T$ is $6$ points, and $W$ is $\CP^2$ blow-up $T$. But $T$ is \emph{not} generic since the $6$ points lie on a conic. Thus, $W$ is \emph{not} a cubic surface. Instead, one can check that $W$ is the resolution of a cubic surface with an ordinary double point. (cf. \cite{reid96} Ch $1$, Ex $13$)}

Now fix some $k > 0$ and assume that $kL + D$ is ample. Let $\bundle{L}^k = \bundle{L}^{\tensor k}$ is the $k$-th tensor power of $\bundle{L}$, and $\bundle{L}_D^k = \bundle{L}^k \tensor \bundle{D} \dsum \bundle{D}$. Let $s''_D$ be a generic section of $\bundle{L}^k \tensor \bundle{D}$, then $s''_D$ defines a generic hypersurface $D'_0$ in the linear system $\linearsys{kL + D}$, and $s_k = \frac{s''_D}{s_D}$ is a (meromorphic) section of $\bundle{L}^k$. Define a ``multi-section'' $r$ of $\bundle{L}$ from $s_k$ as:
$$\Big(r(p)\Big)^{\tensor k} = s_k(p).$$
Consider $r_{\lambda, D} = (\lambda_1 r, \lambda_0)$. Let $W_{\lambda}^k$ be the closure of the image of $r_{\lambda, D}$ in $E$, then similar to \citeequ{inter}, we have
\equlabel{interk}{W_{\lambda}^k \inter V_0 = D'_0, ~~ W_{\lambda}^k \inter V_{\infty} = D_{\infty}, ~~\textrm{ and } W_{\lambda}^k \cdot f = k.}
This is very similar to the ``Branched Covering Trick'' (cf. \cite{bpv}) in which one constructs a branched covering of $X$ with branching loci $B$ in $\linearsys{kL}$ for some divisor $L$. Here, $W_{\lambda}^k$ is generically a branched $k$-ple covering of $X$, and contains the whole fiber $\CP^1$ over the intersection $D'_0 \inter D_{\infty}$. In fact, $W_{\lambda}^k$ can be defined by a generic section of the line bundle $\sheafo_E(k) \tensor \pb{\pi}\bundle{D}$ over $E$. For $k = 1$, this is the case of \citelem{blowup}. (cf. \cite{hart}\S II.7)

\subsection{Semi-stable degeneration of algebraic surfaces}

Let $\polytp = \{(x_1, x_2, x_3) | x_i \geq 0$, for $i = 1, \ldots, 3$, and $\sum_{i = 1}^3 x_i \leq d\} \subset \R^3$, then it defines $\CP^3$ with the divisor $dH$. A generic section of $\sheafo(d)$ defines a generic degree $d$ surface $S$ in $\CP^3$. Every non-singular partition of $\polytp$ give raise to a semi-stable degeneration of $S$. We'll describe some obviously ``universal'' ways to perform such partitions, and the corresponding semi-stable degenerations. Then we'll treat specifically the degeneration of $K3$ surfaces, since it leads to some interesting generalizations in higher dimensional Calabi-Yau hypersurfaces.

\subsubsection{Degeneration of degree $d$ surface into a chain of rational surfaces}

The first kind of partitions are defined by $d-1$ hyperplanes in $\R^3$:
$$L_j = \{(x_1, x_2, x_3) | \sum_{i = 1}^3 x_i = j\}, ~\textrm{ for }~ j = 1, \ldots, d-1.$$
We'll use $\Gamma_d$ to denote this kind of partitions. Label the subpolytope between $L_{j-1}$ and $L_j$ as $\polytp_j$ for $j = 1, \ldots, d$. Let $Y_j$ be the toric variety defined by $\polytp_j$, then, except that $Y_1 \iso \CP^3$, we have $Y_j \iso \blowup{\CP^3}$, the blowup of $\CP^3$ at one point, for $j = 2, \ldots, d$. Let $H$ be the pullback of the hyperplane class, and $E$ be class of the exceptional divisor, then they generate $\coho{2}{\blowup{\CP^3}}$. Let $D_j$ be the divisor on $Y_j$ defined by $\polytp_j$. Then
$$D_1 = H, ~~ \textrm{ and } ~~ D_j = jH - (j-1)E ~ \textrm{ for } j \neq 1.$$
It's well known that $\blowup{\CP^3} \iso \p_{\CP^2}(\sheafo(-1) \dsum \triv)$ is a $\CP^1$-bundle over $\CP^2$. Let $f$ denote the class of a fiber $\CP^1$, then we have the intersections
\equlabel{degDdim2inter}{H \cdot f = E \cdot f = 1 \gives D_j \cdot f = j - (j - 1) = 1.}
Corresponding to the $d$ subpolytopes, there are $d$ components $S_j, ~j = 1, \ldots, d$ in the central fiber of the degeneration of $S$. Each $S_j$ is a toric hypersurface of the corresponding $Y_j$ defined by $D_j$. The intersections $D_j \cdot H$ and $D_j \cdot E$ in $H \iso E \iso \CP^2$ are curves of degree $j$ and $j-1$ respectively. Thus, by \citelem{blowup}, $X_j$ is the blow-up of $\CP^2$ at the intersection of degree $j$ and $j - 1$ curves, which is $j(j - 1)$ points. It's clear from the above description that the intersection of the components $X_j$ and $X_{j+1}$ is a degree $j$ curve, whose self-intersection number is $j$ in $X_j$, for all $0 < j < d$.

The partitions $\Gamma_d$ are special cases of the following. Let $\Gamma_d^k$ be partition defined by the hyperplanes $L_j^k$ in $\R^3$:
$$L_j^k = \{(x_1, x_2, x_3) | \sum_{i = 1}^k x_i = j\}, ~\textrm{ for }~ j = 1, \ldots, d-1,$$
for some fixed $1 \leq k \leq 3$. Notice that $k = 1$ and $k = 3$ give isomorphic partitions in the sense that one is transformed to another by an integral affine isomorphism that maps $\polytp_d$ to itself. Thus $k = 2$ is the only new case. The following pictures show the two partitions for $d = 4$:

\begin{center}
\begin{tabular}{c c c}
\hspace{4.0cm} & \hspace{1.0in} & \hspace{4.0cm} \\
\draw{
 \linewd 0.05
 \rlvec (4 0) \rlvec (-4 4) \rlvec (0 -4)
 \linewd 0.01
 \rlvec (1.732 1) \rlvec (2.268 -1) \move (1.732 1) \rlvec (-1.732 3)

 \linewd 0.03 \move (1 0) \lvec (0 1)
 \linewd 0.02 \lvec (0.433 0.25) \lvec (1 0)

 \linewd 0.03 \move (2 0) \lvec (0 2)
 \linewd 0.02 \lvec (0.866 0.5) \lvec (2 0)

 \linewd 0.03 \move (3 0) \lvec (0 3)
 \linewd 0.02 \lvec (1.299 0.75) \lvec (3 0)
} &
 &
\draw{
 \linewd 0.05
 \rlvec (4 0) \rlvec (-4 4) \rlvec (0 -4)
 \linewd 0.01
 \rlvec (1.732 1) \rlvec (2.268 -1) \move (1.732 1) \rlvec (-1.732 3)

 \linewd 0.03 \move (1 0) \lvec (0 1)
 \linewd 0.02 \lvec (1.299 1.75) \lvec (2.299 0.75) \lvec (1 0)

 \linewd 0.03 \move (2 0) \lvec (0 2)
 \linewd 0.02 \lvec (0.866 2.5) \lvec (2.866 0.5) \lvec (2 0)

 \linewd 0.03 \move (3 0) \lvec (0 3)
 \linewd 0.02 \lvec (0.433 3.25) \lvec (3.433 0.25) \lvec (3 0)

} \\
$\Gamma_4$ & & $\Gamma_4^2$ \\
\end{tabular}
\end{center}

We label the polytope between $L_{j-1}^2$ and $L_j^2$ as $\polytp_j^2$, and let $Y_j^2$ be the toric variety defined by $\polytp_j^2$, for $1 \leq j \leq d$. Then $Y_j^2$ is isomorphic to $\CP^3$ blow-up two non-intersecting line for $1 < j < d$, and is isomorphic to $\CP^3$ blow-up a line for $j = 1, d$. An equivalent description of the $Y_j^2$'s is the following:
$$Y_j^2 \iso Y^2 = \p_{\CP^1 \cross \CP^1}(\pb{p_1}\sheafo(1) \tensor \pb{p_2}\sheafo(-1) \dsum \sheafo), ~\textrm{ for }~ 1 < j < d, ~\textrm{ and }$$
$$Y_1^2 \iso Y_d^2 \iso Z = \p_{\CP^1}(\sheafo(1) \dsum \sheafo(1) \dsum \sheafo),$$
where $p_1$ ($p_2$) are projections to the first (resp. second) factor.

First consider $1 < j < d$, and we denote the divisor class of the $0$-section ($\infty$-section) of the $\CP^1$ bundle as $E$ (resp. $F$), then the divisor defined by $\polytp_j^2$ is:
$$D_j^2 = j E - (j-1) F + d H' = (d+j-1) F - (d-j) E + d H'',$$
where $H = H' + E = H'' + F$ is the pullback of the hyperplane class in $\CP^3$ via the blowdown map. Then the component $X_j^2$ of the central fiber of the degeneration is a toric hypersurface of $Y_j^2$ in the linear system $\linearsys{D_j^2}$. Let $f$ be the class of the fiber then we have the following:
$$E \cdot f = F \cdot f = 1, H' \cdot f = 0 \gives D_j^k \cdot f = 1,$$
which shows that $X_j^2$ is generically a section. Thus by \citelem{blowup}, $X_j^2$ is the blowup of $\CP^1 \cross \CP^1$ at some points. Note that $E \cdot F = 0$, we have $F \cdot F = (H' - H'') \cdot F$, and $E \cdot E = (H'' - H') \cdot E$, thus
$$D_j^2 \cdot E = j H'' \cdot E + (d-j) H' \cdot E, ~~~~~~ D_j^2 \cdot F = (j-1) H'' \cdot F + (d-j+1) H' \cdot F.$$
Let $H_1 = \pb{s} H''$, and $H_2 = \pb{s} H'$, where $s$ stands for $0$ or $\infty$-section of the $\CP^1$ bundle, then $H_i$ is the divisor class of a fiber of $p_i$ in $\CP^1 \cross \CP^1$.  Now \citelem{blowup} tells us the points to be blown-up are the intersection:
$$(j H_1 + (d-j)H_2) \cdot ((j-1) H_1 + (d-j+1) H_2) = j(d-j+1) + (d-j)(j-1) ~\textrm{ points }.$$
Let $C_j$ be a curve in the linear system $\linearsys{j H_1 + (d-j)H_2}$ on $\CP^1 \cross \CP^1$, then $C_j$ has genus $g_j = (d-j-1)(j-1)$ by genus formula. Thus the components $X_j^2$ and $X_{j+1}^2$ intersect along a curve of genus $g_j$.

For $j = 1$, $\polytp_1^2$ is shown below, with $S_i$'s labelling the facets.
\includeps{width=0.25\textwidth}{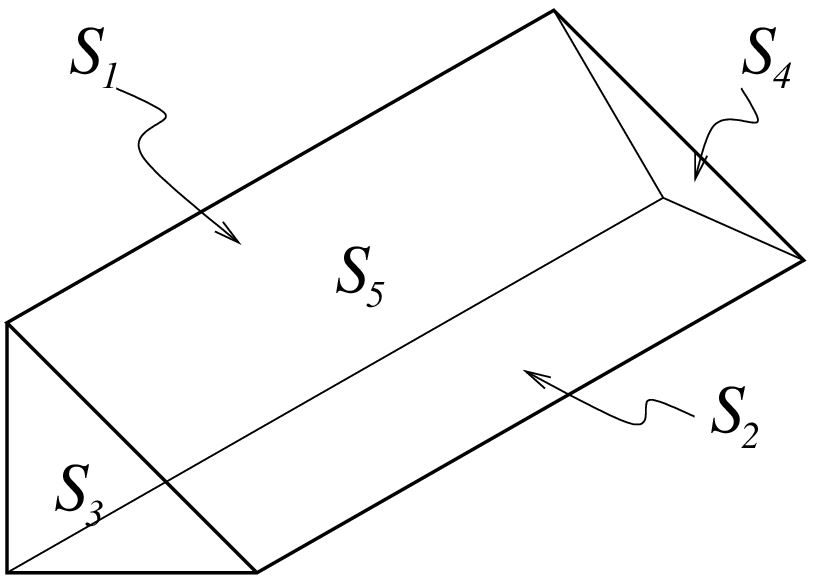}
We'll also use $S_i$'s to denote the divisors corresponding to the faces. With the above notations, we have $E = S_5$, and $H = E + H' = S_5 + S_4 = S_1$ is the pullback of the hyperplane class from $\CP^3$, and the divisor defined by $\polytp_1^2$ is:
$$D_1^2 = E + d H' = S_5 + d S_4.$$
Then $X_1^2$ is a hypersurface in the linear system $\linearsys{D_1^2}$.
\prop{ends}{ $X_1^2$ is in the list $\{\f_d, \f_{d-2}, \ldots \}$.

\proof{
Let $\ell_1 = S_4 \cdot S_5$, $\ell_2 = S_5 \cdot S_1$. Then by simple calculation, we see that 
$$D_1^2 \cdot S_4 = \ell_1, ~~~D_1^2 \cdot S_5 = (d-1) \ell_1 + \ell_2 ~~\textrm{ and } ~~ D_1^2 \cdot S_1 = d \ell_1 + \ell_2.$$
The first equation tells us that $X_1^2$ is a $\CP^1$ bundle over $\CP^1$. The last two intersections are the sections of this $\CP^1$ bundle. Calculating the intersection numbers of the two sections, we see that $X_1^2$ is a Hirtzebruch surface in the list $\{\f_{d+2i}, ~~ i \in \Z\}$. We need more information to pin down the $i$, which can be seen in the following way of determining $X_1^2$.

Consider the divisor $D_1^2 = E + d H'$ in $Y^2$, then a generic hypersurface $\blowup{X_1^2}$ in $\linearsys{D_1^2}$ is generically a section. By \citelem{blowup}, we see that $\blowup{X_1^2}$ is a blowup of $\CP^1 \cross \CP^1$ at $(H_1 + (d - 1)H_2) \cdot d H_2 = d$ points. Let $\ell = \numin{\polytp_1^2}$, then morphism $\phi_1^2$ of $Y^2$ to $\CP^{\ell - 1}$ defined by $\polytp_1^2$ contracts the fiber of $p_2$ in $F$, from which $Y_1^2$ appears as the image. The image of generic $\blowup{X_1^2}$ is a (smooth) hypersurface of $Y_1^2$, which is $X_1^2$. Since the map $Y^2 \to Y_1^2$ is a blowup, we see that $\blowup{X_1^2} \to X_1^2$ is a blow-up as well. In fact, since $D_1^2 \cdot F = d H_2$, $\blowup{X_1^2}$ is the blowup of $X_1^2$ at $d$ points. Thus, $X_1^2$ is rational, and fits into the blowing-up diagram
$$\CP^1 \cross \CP^1 \leftarrow \blowup{X_1^2} \to X_1^2.$$
In general, the blowing-up blowing-down process as above can only tell us that $X_1^2$ is in the list $\{\f_d, \f_{d-2}, \ldots\}$.
}
}

Similar arguments as above apply for $j = d$. Thus, the partition $\Gamma^2$ gives rise to a semi-stable degeneration of degree $d$ surface into a chain of rational surfaces, with both ends in the list $\{\f_d, \f_{d-2}, \ldots\}$.

\subsubsection{Inductive structure of degree $d$ surfaces}

We can partition the polytope $\polytp_d$ using some of the hyperplanes $L_j$ as above, and get a partition which has fewer subpolytopes. One sepcial case is to take the partition of $\polytp_d$ given by a single hyperplane $L_{d-1}$. We end up with a semi-stable degeneration of degree $d$ surface into two components. One of the component is a degree $d-1$ surface, and the other is the $X_d$ described above, which is the blow-up of $\CP^2$ at $d(d-1)$ points. Thus, surfaces of all degrees in $\CP^3$ as a system has an ``inductive'' structure.

Of cause, we can partition $\polytp_d$ using some of the hyperplanes $L_j^2$ as well. The next sub-section provides some interesting examples in the case of $K3$ surfaces.

\subsubsection{Degeneration of $K3$ surfaces}

First, we partition the polytope $\polytp_4$ into $2$ pieces by $L_2^2$, and denote the corresponding toric varieties as $Z_1$ and $Z_2$. Then $Z_1 \iso Z_2 \iso Z$. The subpolytopes define the same divisor in $Z$:
$$D = 2 S_5 + 4 S_4.$$
This partition induces a semi-stable degeneration of $K3$ surface into $2$ components $P_i$, $i = 1, ~ 2$. Then $P_i$ are in the same linear system $\linearsys{D}$, thus $P_1 \iso P_2 \iso P$.
\prop{K3toE1}{ $P \iso E(1)$ as topological manifolds.

\proof {
It's easy to compute the intersection $D \cdot S_5 = 2 \ell_1 + 2 \ell_2$, then by genus formula we see that $P$ intersects with $S_5$ in an elliptic curve $A$. By adjunction
$$K_P = (K_{Z} + D) \restricto{P} = (- S_5)\restricto{P},$$
we see that $K_P = - A$, and
$$c_1^2(P) = A \cdot A = (S_5)^2 \restricto{P} = (S_5)^2 \cdot D = 0.$$
By Lefschetz hypersurface theorem, we have $q(P) = 0$. Since $-K_P = A$ is effective, $K_P \cdot K_P = 0$, and $q(P) = 0$, we see that $P$ is rational by the Castelnuovo-Enriques theorem.
Using the cohomology sequence associated to the exact sequence
$$\shortexact{0}{\sheafo_{Z}(-D)}{}{\sheafo_{Z}}{}{\sheafo_P}{0}$$
and $h^2(\sheafo_{Z}) = h^3(\sheafo_{Z}) = 0$, we see that 
$$h^2(\sheafo_P) = h^3(\sheafo_{Z}(-D)) = 0,$$
while the latter is $0$ by Serre duality and Kodaira vanishing. By Noether's formula, we get $e(P) = c_2(P) = 12$, and thus $h^2(P) = h^{1,1}(P) = 10$. Then the signature
$$\sigma(P) = \frac{1}{3}(c_1^2(P) - 2 c_2(P)) = -8.$$
Let $Q$ be the intersection form of $P$, then $Q = \langle 1 \rangle \dsum 9 \langle -1 \rangle$ if it's odd, or $Q = -E_8 \dsum H$ if it's even where $H$ is the intersection form of $2$-torus. From 4-manifold theory, $Q \neq -E_8 \dsum H$, which gives
$$Q = \langle 1 \rangle \dsum 9 \langle -1 \rangle.$$
We see that $P$ is the blow-up of some ruled surface $\f_q$ at $8$ points. Thus, topologically, we have $P \iso E(1)$.
}
}
In light of \citeprop{K3toE1}, we may think of the degeneration as an algebraic geometrical version of the fiber sum operation.

Next, we consider the partition by $L_3^2$. There are two components in the central fiber of the induced semi-stable degeneration of $K3$ surfaces. Both components are hypersurfaces in $Z$. One of them, say $W_1$, is in the linear system $\linearsys{S_5 + 4 S_4}$, which by \citeprop{ends} is one of $\f_4, \f_2$, and $\CP^1 \cross \CP^1$. The other, $W_2$, is in the linear system $\linearsys{3 S_5 + 4 S_4}$. By toric geometry, the canonical class of $Z$ is
$$K_Z = - \sum_{i = 1}^5 S_5 = -(3 S_5 + 4 S_4).$$
Thus, $K_{W_2} = (W_2 + K_Z) \restricto{W_2} = 0$ by adjunction, and $h^1(W_2) = 0$ by Lefschetz hypersurface theorem. We see that $W_2$ is again a $K3$ surface. The two components intersect along a rational curve $\ell$ by the calculation in the first subsection, and the self-intersection of $\ell$ in $W_2$ is $-2$ either by checking directly or using genus formula.

In next subsection, we'll construct the classic degeneration of $K3$ into $4$ cubic surfaces, and the generalization of it together with the first construction above to higher dimensional Calabi-Yau hypersurfaces.

\subsection{Semi-stable degeneration of Calabi-Yau hypersurfaces}
\subsubsection{Classical degeneration of $K3$ and its generalizations}
Let $\polytp \subset \R^n$ (resp. $\Gamma$) be the polytope (resp. partition) as given in \citexmp{degenCPn}. We'll adopt the notations of \citexmp{degenCPn}.
The $Q^n$ in \citexmp{degenCPn} is defined by the polytope $\polytp_0(n) \subset \R^n$. The $1$-skeleton of the dual fan $\fan_0(n)$ of $\polytp_0(n)$ is generated by $2n$ vectors: $e_1, \ldots, e_n$ and 
\equlabel{normals}{w_j = -\sum_{i = 1}^j e_i, ~ \textrm{ for } ~ j = 1, \ldots, n.}
There are $2^n$ $n$-dimensional cones in $\fan_0(n)$, and they are generated by one of $\pm e_1$ with certain $n-1$ of the other $1$-dimensional cones in the above list. Let $D_1, \ldots, D_n$ be the divisors corresponding to $e_1, \ldots, e_n$, and $D_{n+1}, \ldots, D_{2n}$ be the divisors corresponding to $w_1, \ldots, w_n$, respectively. Then $\polytp_0(n)$ defines divisor
$$D = \sum_{i = 1}^n D_i$$
on $Q^n$. Sometimes we will use $D_i(n)$ to emphasize that $D_i$ is a divisor of $Q^n$. We also use $H(n)$ to denote $D_1(n)$. The properties of $Q^n$ are collected in the following:

\prop{aboutQn}{Fix any $n > 0$, and set $Q^0 = pt$, we have the following:
\begin{enumerate}
\item $Q^n$ is rational.
\item $Q^n$ is Fano, with canonical class $K = -\sum_{i = 1}^{2n} D_i$.
\item $D_j = \sum_{i = j}^n D_{n+i}$ for $j = 1, \ldots, n$.
\item $Q^{n} = \p_{Q^{n-1}}([H(n) \dsum \triv])$.
\item $D_1$ and $D_{n+1}$ are sections of the above bundle.
\item $D_1 \restricto{D_1} = H(n-1)$, and $D_{n+1} \restricto{D_{n+1}} = -H(n-1)$,
\item Let $\pi_n : Q^{n} \to Q^{n-1}$ be the bundle map, then
$D_{i+1}(n) = \inverse{\pi_n}(D_i(n-1))$, for $i = 1, \ldots, n-1$, and
$D_{i+2}(n) = \inverse{\pi_n}(D_i(n-1))$, for $i = n, \ldots, 2(n-1)$.
\item $\hodgenum{p}{p} = \combn{n}{p}$ for all $p$.
\end{enumerate}

\proof{
Everything is quite straightforward using direct toric arguments.
}
}

Now we should describe the degeneration of a Calabi-Yau hypersurface $V$ in $\CP^n$. The components in the central fiber of the degeneration of $V$ will consists of $n+1$ hypersurfaces of the Fano variety $Q^n$ in the linear system $\linearsys{\sum_{i = 1}^n D_i}$. It's easy to see from the partition $\Gamma$ that the dual graph of the central fiber ($V_0$) is a triangulation of $S^{n-1}$ by $n+1$ simplices. From direct toric calculation, we see that $D = \sum_{i = 1}^n D_i$ is generically a section of the $\CP^1$ bundle $Q^{n} \to Q^{n-1}$. Apply \citelem{blowup}, we see that the components are all rational. To know what exactly those components are, we need to find the intersections of $D$ with the $0$- and $\infty$-sections.

Let's first consider what happens for $K3$ surfaces in $\CP^3$. Recall $Q^3 = \p_{Q^2}([H] \dsum \triv)$. The class $D = D_1 + D_2 + D_3$, where $D_1$ is the class of $0$-section. Then
$$D \cdot D_1 = (D_1 + D_2 + D_3) \cdot D_1 = H + H + f = 3H - E$$
as divisor in $Q^2 = \blowup{\CP^2}$. Note that $3H - E = H + 2f + E$ is the anticanonical divisor of $Q^2$. The class of $\infty$-section is $D_4$, and 
$$D \cdot D_4 = (D_1 + D_2 + D_3) \cdot D_4 = H + f = 2H - E$$
as divisor in $Q^2$. Then the class of $T$ as in \citelem{blowup} is
$$(2H - E) \cdot (3H - E) = 5 ~\textrm{ points}.$$
Thus, the components are $Q^2$ blow-up $5$ points, which is $\CP^2$ blow-up $6$ points. By contracting $E$, we see that the $6$ points in $\CP^2$ are the intersection of two curves of degree $2$ and $3$. As in \citermk{blowup1}, we get the resolution of a cubic hypersurface with an ordinary double point.

Since the $Q^3$'s in the central fiber of the degeneration of $\CP^3$ meet along the subvarieties $D_j$ $j = 4, 5, 6$, we can get the singular loci of $V_0$ by computing the intersections of the components with $D_j$ for $j = 4, 5, 6$. For $D \cdot D_4 = 2H - E$, we see that the intersection is a $\CP^1$. $D_5$ is also isomorphic to $\blowup{\CP^2}$ and $D \cdot D_5 = H + f = 2H - E$, gives a $\CP^1$ as well. As to $D^6$, it's isomorphic to $\CP^1 \cross \CP^1$, and $D \cdot D_6 = P_1 + P_2$, where $P_i$ is the class of the $i$-th factor ($i = 1, 2$), and we get another $\CP^1$. Thus, the components of $V_0$ intersect along $\CP^1$'s, and we recover the known construction of $K3$ degenerates into $4$ rational surfaces.

Next, let's consider the case of quintics in $\CP^4$, and the general case would not be much different. The components of $V_0$ are hypersurfaces of $Q^4$ in the linear system $\linearsys{D = \sum_{i = 1}^4 D_i}$. In the $\CP^1$-bundle $Q^4 = \p_{Q^3}([H(3) \dsum \triv]) \to Q^3$, $D_1$ is the $0$-section, and $D_5$ is the $\infty$-section. Using \citeprop{aboutQn} and let $K(3)$ be the canonical class of $Q^3$. We have the following intersections:
$$D \cdot D_1 = (\sum_{i = 2}^8 D_i) \cdot D_1 = -K(3),$$
$$D \cdot D_5 = -K(3) - H(3) = \sum_{i = 1}^3 D_i(3),$$
where $D_i(3)$'s are classes in $Q^3$. Thus, the $T$ as in \citelem{blowup} is the intersection of a $K3$ surface and a cubic surface in $Q^3$. One can check that the curve $T$ has genus $12$, and is in the class $8H(2) -3E(2) + 3f$ where $H(2)$ (resp. $E(2)$) is the image of the class $H$ (resp. $E$) in $Q^2$ under the embedding $Q^2 = D_1(3) \into Q^3 = D_1(4) \into Q^4$, while $f$ is the class of fiber in $Q^4$.

Let $C^n$ denote one component in the degeneration of $n$-dimensional Calabi-Yau hypersurface in $\CP^{n+1}$, then $C^n$ is in the linear system $\linearsys{\sum_{i = 1}^{n+1} D_i(n+1)}$ of $Q^{n+1}$. $C^n$ is a blow-up of $Q^{n}$ along the intersection of two hypersurfaces in the linear systems $\linearsys{\sum_{i = 1}^{2n}D_i(n)}$ and $\linearsys{\sum_{i = 1}^n D_i(n)}$ respectively. The first one is a Calabi-Yau since it's in the anticanonical class, and the second one is in fact $C^{n-1}$. We'll compute the class of the intersection in the following. Let
$D_{\bar{1}}(n) = \sum_{i = 2}^n D_i(n),$
then
$$D_{\bar{1}}(n) = \inverse{\pi_{n}}(\sum_{i = 1}^{n-1} D_i(n-1)) = \inverse{\pi_n}(C^{n-1}).$$
Now, since $\sum_{i = 1}^{2n}D_i(n) = 2D_1(n) + D_{\bar{1}}$, we have
$$\equlist{rcl}{
 T & = & \sum_{i = 1}^{2n}D_i(n) \cdot \sum_{i = 1}^n D_i(n) \\
 & = & 2 D_1(n) \cdot D_1(n) + 3 D_1(n) \cdot D_{\bar{1}}(n) + D_{\bar{1}}(n) \cdot D_{\bar{1}}(n) \\
& = & 2 D_1(n-1) + 3 C^{n-1} + \inverse{\pi_n}(C^{n-1} \cdot C^{n-1}).
}$$
where $D_1(n-1)$ and $C^{n-1}$ are regarded as classes of $Q^n$ through the embedding $Q^{n-1} = D_1(n) \into Q^n$. Then we can carry on the computation for $C^{n-1} \cdot C^{n-1}$, and eventually we arrive at an expression involving $D_1(j)$, $C^{j}$, and composition of pull-backs by $\pi_j$.

We summarize what we got as the following:
\thm{CYdegen}{There exists a semi-stable degeneration of degree $n+1$ Calabi-Yau hypersurface in $\CP^n$, such that the dual graph of the central fiber is the triangulation of $S^{n-1}$ by the boundary of $n$-simplex. The components in the central fiber are all isomorphic to a blow-up of $Q^{n-1}$ along the intersection of two hypersurfaces.}

\subsubsection{Degeneration of Calabi-Yau hypersurfaces in $\CP^{2k+1}$}

In this subsection, we generalize the degeneration of $K3$ surface into two copies of $E(1)$. Although it's hard to determine the components in the central fiber, we find that the many properties generalize.

Calabi-Yau hypersurfaces in $\CP^{2k+1}$ are in the linear system defined by the polytope $\polytp_{2k+2}$. Consider the partition defined by $L_{k+1}^{k+1}$. Then there are two subpolytopes $\polytp_1$ and $\polytp_2$, and we denote the corresponding toric varieties as $X_1$ and $X_2$. It's easy to check that
$$X_1 \iso X_2 \iso X = ~\textrm{ the blowup of}~ \CP^{2k+1} ~\textrm{ along linear subspace }~ \CP^{k}.$$
Equivalently, we have $X = \p_{\CP^k}(\sheafo^{k+1}(1) \dsum \triv)$. The exceptional divisor is $E = \CP^k \cross \CP^k$, and the canonical class of $E$ is $K_E = - (k+1) (\pb{p_1} H_1 + \pb{p_2} H_2)$, where $p_i$ and $H_i$ are the projections to and the hyperplane class on the $i$-th factor respectively ($i = 1, 2$). Let $H = H' + E$ be the pull back of the hyperplane class in $\CP^{2k+1}$, then simple calculation gives:
$$K_X = -(2k+2) H' - (k+2) E, ~\textrm{ and }~ E \cdot E = (H - H') \cdot E = \pb{p_1} H_1 - \pb{p_2} H_2,$$
where we view $\pb{p_i} H_i$ as classes in $X$ via the inclusion $E \into X$.

The induced semi-stable degeneration of Calabi-Yau hypersurface has two components in the central fiber, say $Y_i \subset X_i$ ($i = 1, 2$). By toric computation, the divisors defined by $\polytp_i$ ($i = 1, 2$) are linearly equivalent to $D = (k+1)(H' + H)$, thus $Y_1 \iso Y_2 \iso Y$, and the canonical class of $Y$ is $K_Y = (-E) \restricto{Y}$. By $D\restricto{E} = (k+1)(H' + H) \restricto{E} = (k+1)(\pb{p_1}H_1 + \pb{p_2}H_2) = -K_E$, we see that the anticanonical class $-K_Y$ defines a $2k-1$ dimensional Calabi-Yau $Z$ in $Y$. The two components $Y_i$ intersect along $Z$. By property of semi-stable degeneration, the self-intersection of $Z$ in the components add up to $0$, while $Z$ has same self-intersection in this case, we see that $Z \cdot Z = 0$, i.e. $(K_Y)^2 = 0$. We can state the result as:

\thm{degenCYodd}{There exists a semi-stable degeneration of Calabi-Yau hypersurface in $\CP^{2k+1}$ into two components $Y_i$ ($i = 1, 2$), such that $Y_1 \iso Y_2 \iso Y$ satisfy the following: $(K_Y)^2 = 0$, $-K_Y$ is effective and defines a Calabi-Yau subvariety $Z$ of $Y$, $Z$ is isomorphic to a hypersurface in $\CP^k \cross \CP^k$, and the singular loci of the central fiber is $Z$.}

\subsubsection{Degeneration of Calabi-Yau hypersurface in $\p(1,1,2,2,2)$}

By \citethm{weakdegen} and \citexmp{degenWPn}, we have a similar weak semi-stable degeneration of the Calabi-Yau hypersurface in $\p(1,1,2,2,2)$ into $5$ components. There is one smooth component which is isomorphic to $C^{4}$ above. The other four components are toric hypersurfaces with quotient singularities.

\subsection{Semi-stable degeneration of higher dimensional toric hypersurfaces}

The construction described in $5.3$ can be easily generalized to hypersurfaces in $\CP^n$ with $n > 3$. The partition $\Gamma_d^k$ for $1 < k < \frac{n + 1}{2}$ is given by the hyperplanes:
$$L_j^k = \{(x_1, \ldots, x_n) | \sum_{i = 1}^k x_i = j\}, ~\textrm{ for }~ j = 1, \ldots, d-1.$$
We'll briefly describe the special case $k = 1$. Similar to the surface case, the partition $\Gamma_d$ induces a semi-stable degeneration of degree $d$ hypersurface $X$ into a chain of rational varieties $X_i$, $i = 1, \ldots, d$. For each $i$, $X_i$ is the blowup of $\CP^{n-1}$ along the intersection of two hypersurfaces of degree $i$ and $i-1$ respectively. If we define a partition using only $L_{d-1}$, then we see that degree $d$ hypersurfaces in $\CP^n$ as a system has some ``inductive'' structure.

\section{Some generalizations}

We describe a construction of liftings with higher extra dimensions. One example of such construction recovers the special standard model of expanded degeneration in \cite{lijun}.

\subsection{Semi-stable degeneration with higher dimensional base}
We can construct higher dimensional liftings for certain partitions, which lead to degernation with higher dimensonal base. For simplicity, we'll assume that all polytopes are non-singular, and partitions are non-singular as well.

Let $\polytp$ be a non-singular polytope. Suppose partition $\Gamma$ is defined by parallel hyperplanes $L_i = \{\pairing{m}{n} = c_i\}$, where $n \in N$ is a fixed primary vector, and $c_1 < c_2 \ldots, < c_l$ are integers. Let $\{\polytp_i\}_{i = 0}^{l}$ be the subpolytopes in $\Gamma$, such that $\polytp_i$ lies between $L_{i}$ and $L_{i+1}$ when both are defined. We perform the lifting in $l$ steps as following.

Step $1$: Lift $\polytp$ by the partition defined by $L_1$, to $\polytp_1 \subset \polytp \cross \R \subset M_{\R} \cross \R$. Let $L_{i; 1}$ be the preimage of $L_i$ ($i > 1$) under the projection $\pi_1 : M_{\R} \cross \R \to M_{\R}$, then $L_{i; 1}$ ($i > 1$) define a partition $\Gamma_1$ of $\polytp_1$. It's easy to check that $\Gamma_1$ is non-singular, given everything we started with is non-singular.

Step $2$: Repeat Step $1$ for $\polytp_1$, $\Gamma_1$, and $L_{2; 1}$ to get $\polytp_2 \subset M_{\R} \cross \R^2$, and partition $\Gamma_2$ of $\polytp_2$ given by $(l-2)$ hyperplanes $L_{i; 2}$ ($i > 2$) which are the preimage of $L_{i}$ ($i > 2$) under the projection $\pi_2 : M_{\R} \cross \R^2 \to M_{\R}$.

$\ldots$

Step $l$: The result is a non-singular polytope $\polytp_l \subset M_{\R} \cross \R^l$.\\
In each step above, the lifting polytope $\polytp_i$ can be chosen to be either compact or open. The following are the results of the first $2$ steps for the partition on $[0, 4]$ defined by lattice points $1, 2, 3$.
\includeps{width=0.75\textwidth}{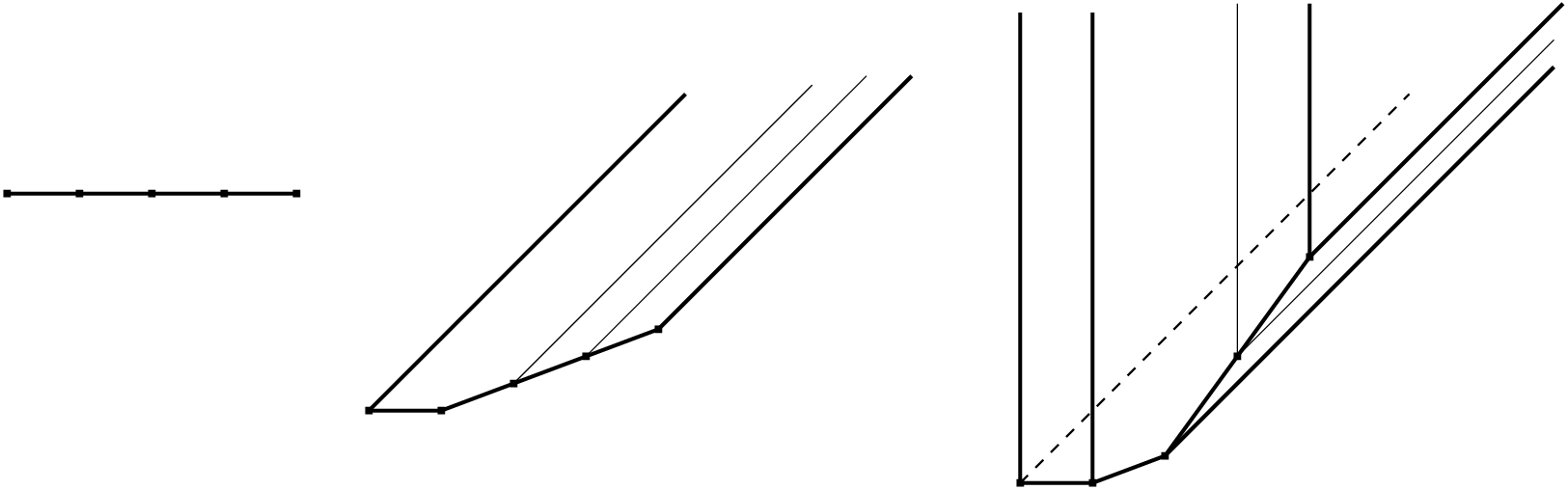}

Another way to get the final result in one step is to construct a suitable function $F_{\Gamma, l}: \polytp \to \R^l$, and let $G_{\Gamma, l}$ be the graph of $F_{\Gamma, l}$ in $\polytp \cross \R^l$. Then we get polytope $\polytp_l$ by taking the Minkowski sum of $G_{\Gamma, l}$ with $(\Rplus)^l$. The function $F_{\Gamma, l}$ can be constructed in steps that are parallel to the steps above. Instead, we'll define it as following. Let $\{e_i\}|_{i = 1}^l$ be the standard basis of $\R^l$. For each $L_i$, define $f_i = \pairing{m}{n} - c_i$, then $f_i \restricto{L_i} = 0$. Let $f_0 = 0$. Define
$$F_{\Gamma, l} = \sum_{j = 0}^{i} f_i e_i, ~\textrm{ on }~ \polytp_i ~\textrm{ for }~ i = 0, \ldots, l,$$
then it's easy to check that $F_{\Gamma, l}$ induces the polytope $\polytp_l$.

Let $X$ (resp. $\new{X}_l$) be the toric variety defined by $\polytp$ (resp. $\polytp_l$). By similar argument as in the proof of \citethm{degeneration}, there is a morphism of toric varieties $p_l : \new{X}_l \to \C^l$, with genric fiber $X$, and singular fibers over the coordinate planes in $\C^l$. The singular fibers correspond to partitions defined by subsets of $\{L_i\}_{i = 1}^l$, especially, the fiber over $(0, \ldots, 0)$ corresponds to the original partition $\Gamma$. Again, we want to consider the induced degeneration of subvarieties of $X$. We'll consider compact liftings in the following, and use the basic construction with compact polytopes.

\lemma{convex}{Suppose $\polytp \subset M_{\R}$ is a non-singular polytope, and $\Gamma$ is a partition defined by a single hyperplane. Let $\new{\polytp} \subset M_{\R} \cross \R$ be a compact lifting of $\polytp$ by $\Gamma$. Let $\fan$ and $\new{\fan}$ be the fans corresponding to $\polytp$ and $\new{\polytp}$ respectively, then $\new{\fan}$ has $\fan$ as a subfan. Suppose $\phi$ is an integral convex piecewise linear function on $\fan$, then there exists an integral convex piecewise linear function $\new{\phi}$ on $\new{\fan}$ that restrict to $\phi$ on $\fan$.

\proof{
By construction, $\new{\fan}$ has three more $1$-dimensional cones than $\fan$. We may assume that $-e_1 = (0, \ldots, 0, -1) \in N_{\R} \cross \R$ generates an edge of $\new{\fan}$, and let $v_1$ and $v_2$ generate the two edges in upper half space. Let $\new{\fan}_i$ be the fan generated by $\fan$ and $v_i$ for $i = 1, 2$, then $\support{\fan_i} = \cl{\halfsp^+} = N_{\R} \cross \cl{\R^+}$. Let $\pi_i: N_{\R} \cross \R \to N_{\R}$ be the projection with kernel $\Span(v_i)$, then the function $\phi_i = \pb{\pi_i} (\phi)$ is convex on the fan generated by $\fan$ and $v_i$, for $i = 1, 2$. By definition, $\phi_1(v_1) = \phi_2(v_2) = 0$. Let
$$\new{\phi} = \max_{t \in [0, 1]}(t \phi_1 + (1-t) \phi_2),$$
then the set $\{(x, y) ~|~ y \leq \new{\phi}(x)\} \subset \cl{\halfsp^+} \cross \R$ is the convex hull of the corresponding sets of $\phi_1$ and $\phi_2$. Thus $\new{\phi}$ is convex on $\cl{\halfsp^+}$. Now, pick $0 \gg a \in \Z$, and set $\new{\phi}(-e_1) = a$, and linearly extend to all $\halfsp^-$ on each cone of $\new{\fan}$. Then $\new{\phi}$ is convex on the whole $N_{\R} \cross \R$. It's obvious that $\new{\phi}$ is piecewise linear on $\new{\fan}$, and integrality follows.
}
}
\rmk{feature}{The function $\new{\phi}$ constructed in the above proof has an extra feature, i.e. $\new{\phi} = 0$ on the cone generated by the new edges $v_i$ in $\new{\fan}$.}

Suppose $D$ is a divisor on $X$ which is generated by sections, then its support function $\phi_D$ is a convex integral function on $\fan$. Apply \citelem{convex} on $\phi_D$ inductively to each step of the (compact) lifting of $\polytp$ to $\polytp_l$, we get a function $\phi_l$ which is integral and convex on the fan $\new{\fan}_l$ given by the compact lifting $\polytp_l$. Let $\face \faceof \polytp$ be a facet (i.e. codimension one face), and $\new{\face} \faceof \polytp_l$ be the facet that maps to $\face$ under the projection $N_{\R} \cross \R^l \to N_{\R}$. Suppose $D = \sum_{\face} a_{\face} D_{\face}$ where $D_{\face}$ is the divisor on $X$ corresponding to $\face$. From the \citermk{feature}, we see the the divisor defined by $\phi_l$ is $\new{D}_l = \sum_{\face} a_{\face} D_{\new{\face}}$, where $D_{\new{\face}}$ is the divisor on $\new{X}_l$ defined by $\new{\face}$. Since $\phi_l$ is convex, the linear system $\linearsys{\new{D}_l}$ is generated by sections. Let $V$ be a generic hypersurface of $X$ in the linear system $\linearsys{D}$, and  $\new{V}$ be a generic hypersurface of $\new{X}_l$ in the linear system $\linearsys{\new{D}_l}$, then $\new{V}$ intersects the invariant divisors transversely. Thus, over a neighbourhood of $(0, \ldots, 0) \in \C^n$, we can restrict $p_l$ to $\new{V}$ and get a semi-stable degeneration of $V$ with a $l$-dimensional base.

\subsection{Degeneration of non-compact varieties}

We need not to restrict ourselves to compact varieties. The construction in section $3$ readily works for non-singular (or mildly singular) partitions on non-singular (resp. simplicial) non-compact polytopes. Since polytope by definition is a finite intersection of half-spaces, we can define a fan using the supporting hyperplanes. Then a toric variety can be constructed from the fan. Thus, we can construct a toric variety from a non-compact polytope. The proof of \citethm{degeneration} goes through without alteration, and we have a degeneration of non-compact toric varieties.

\exmp{degenOpCn}{Let $\polytp = M_{\R}$ the whole space. $\polytp$ is defined equivalently as
$$\{m \in M_{\R} ~|~ \pairing{m}{0} \geq 0\}.$$
Thus the fan of $\polytp$ is $\fan = \{ 0 \}$, the trivial fan in $N_{\R}$. The corresponding toric variety is simply the torus $(\C^*)^n$. Let $\Gamma$ be the partition defined by the fan $\fan(n)$ as in \citexmp{degenCPn}, then it's non-singular. The semistable degeneration defined by $\Gamma$ is the familiar map $p: \C^{n+1} \to \C$, where $p(x_0, x_1, \ldots, x_n) = \prod_{i = 0}^n x_i$.
}
\exmp{degenLJ}{Let $\polytp = \R$, then as above, the variety is $\C^*$. Without loss of generality, we might let $\Gamma$ be defined by the integers $\{j \in \Z ~|~ 0 \leq j \leq l\}$. Carry out the construction in the first subsection, we recover the construction of the special standard model of expanded degeneration in \cite{lijun}.}

\bibliographystyle{amsalpha}

\bibli{99}{
\bitem{alexeev}{V.~Alexeev}{Complete Moduli in the Presence of Semiabelian Group Action}{arXiv:math.AG/9905103 v2}
\bitem{ambramovich}{D.~Ambramovich, K.~Karu}{Weak Semistable Reduction In Characteristic 0, Preliminary Version}{arXiv:alg-geom/9707012}
\bitem{bpv}{W.~Barth, A.~Van de Ven, C.A.M.~Peters}{Compact complex surfaces}{Ergebnisse der Mathematik und ihrer Grenzgebiete. 3. Folge, Bd. 4, Springer-Verlag, 1984}
\bitem{cox}{D.~Cox}{Recent Developements in Toric Geometry}{alg-geom/9606016}
\bitem{eh00}{D.~Eisenbud, J.~Harris}{The Geometry of Schemes}{Graduate Texts in Mathematics 197, Springer, 2000}
\bitem{fulton}{W.~Fulton}{Introduction to Toric Varieties}{Annals of Mathematics Studies 131, Princeton University Press, 1993}
\bitem{hart}{R.~Hartshorne}{Algebraic Geometry}{Graduate Texts in Mathematics 52, Springer-Verlag, 1977}
\bitem{ionelparker}{E.-N.~Ionel, T.H.~Parker}{The Symplectic Sum Formula for Gromov-Witten Invariants}{arXiv:math.SG/0010217}
\bitem{kkms}{G.~Kempf, F.~Knudsen, D.~Mumford and B.~Saint-Donat}{Toroidal Embeddings I}{Lecture Notes in Mathematics 339, Springer-Verlag, 1973}
\bitem{liruan}{A.-M.~Li, Y.~Ruan}{Symplectic surgery and Gromov-Witten invariants of Calabi-Yau 3-folds}{Invent. Math. 145 (2001), no. 1, 151--218}
\bitem{lijun}{J.~Li}{A degeneration of stable morphisms and relative stable morphisms}{arXiv:math.AG/0009097}
\bitem{mcduff}{D.~McDuff, D.~Salamon}{Introduction to Symplectic Topology}{Oxford Mathematical Monographs, Oxford Science Publications, 1995}
\bitem{oda}{T.~Oda}{Convex Bodies and Algebraic Geometry: An Introduction to the Theory of Toric Varieties}{Ergbnisse der Mathematik und ihrer Grenzgebiete. 3. Folge, Bd. 15, Springer-Verlag, 1988}
\bitem{reid96}{M.~Reid}{Chapters on Algebraic Surfaces}{arXiv:alge-geom/9602006}
\bitem{spanier}{E.~H.~Spanier}{Algebraic Topology}{McGraw-Hill, 1966}
}
\end{document}

%% file: lattice.tex
 \linewd 0.02 \arrowheadtype t:V \arrowheadsize l:0.3 w:0.07
 \ravec (2 0) \move (0 0) \ravec (0 2) \move (0 0)
 \fcir f:0 r:0.05
 \rmove (0.5 0) \fcir f:0 r:0.05
 \rmove (0.5 0) \fcir f:0 r:0.05
 \rmove (0.5 0) \fcir f:0 r:0.05
 \move (0 0.5)
 \fcir f:0 r:0.05
 \rmove (0.5 0) \fcir f:0 r:0.05
 \rmove (0.5 0) \fcir f:0 r:0.05
 \rmove (0.5 0) \fcir f:0 r:0.05
 \move (0 1)
 \fcir f:0 r:0.05
 \rmove (0.5 0) \fcir f:0 r:0.05
 \rmove (0.5 0) \fcir f:0 r:0.05
 \rmove (0.5 0) \fcir f:0 r:0.05
 \move (0 1.5)
 \fcir f:0 r:0.05
 \rmove (0.5 0) \fcir f:0 r:0.05
 \rmove (0.5 0) \fcir f:0 r:0.05
 \rmove (0.5 0) \fcir f:0 r:0.05

%% file: latticebg.tex
 \linewd 0.02 \arrowheadtype t:V \arrowheadsize l:0.3 w:0.07
 \ravec (2.5 0) \move (0 0) \ravec (0 2.5) \move (0 0)
 \fcir f:0 r:0.05
 \rmove (0.5 0) \fcir f:0 r:0.05
 \rmove (0.5 0) \fcir f:0 r:0.05
 \rmove (0.5 0) \fcir f:0 r:0.05
 \rmove (0.5 0) \fcir f:0 r:0.05
 \move (0 0.5)
 \fcir f:0 r:0.05
 \rmove (0.5 0) \fcir f:0 r:0.05
 \rmove (0.5 0) \fcir f:0 r:0.05
 \rmove (0.5 0) \fcir f:0 r:0.05
 \rmove (0.5 0) \fcir f:0 r:0.05
 \move (0 1)
 \fcir f:0 r:0.05
 \rmove (0.5 0) \fcir f:0 r:0.05
 \rmove (0.5 0) \fcir f:0 r:0.05
 \rmove (0.5 0) \fcir f:0 r:0.05
 \rmove (0.5 0) \fcir f:0 r:0.05
 \move (0 1.5)
 \fcir f:0 r:0.05
 \rmove (0.5 0) \fcir f:0 r:0.05
 \rmove (0.5 0) \fcir f:0 r:0.05
 \rmove (0.5 0) \fcir f:0 r:0.05
 \rmove (0.5 0) \fcir f:0 r:0.05
 \move (0 2)
 \fcir f:0 r:0.05
 \rmove (0.5 0) \fcir f:0 r:0.05
 \rmove (0.5 0) \fcir f:0 r:0.05
 \rmove (0.5 0) \fcir f:0 r:0.05
 \rmove (0.5 0) \fcir f:0 r:0.05